%% file: milnorfiberflatness.tex
\documentclass[12pt]{article}
\usepackage[margin=1in]{geometry}

\usepackage[style = alphabetic, locallabelwidth]{biblatex}
\usepackage{hyperref}
\usepackage{enumitem}

\usepackage{tikz}
\usepackage{tikz-cd}
\usetikzlibrary {decorations.pathmorphing, decorations.pathreplacing, decorations.shapes}

\usepackage{amsmath}
\usepackage{amssymb}
\usepackage{amsthm}
\usepackage{mathrsfs}

\newcommand{\bC}{\mathbb{C}}
\newcommand{\bE}{\mathbb{E}}
\newcommand{\bF}{\mathbb{F}}

\newcommand{\bP}{\mathbb{P}}
\newcommand{\bQ}{\mathbb{Q}}

\newcommand{\bZ}{\mathbb{Z}}

\newcommand{\cF}{\mathcal{F}}
\newcommand{\cI}{\mathcal{I}}
\newcommand{\cJ}{\mathcal{J}}
\newcommand{\cO}{\mathcal{O}}

\newcommand{\cU}{\mathcal{U}}

\newcommand{\mFm}{\mathfrak{m}}

\DeclareMathOperator{\coker}{coker}
\DeclareMathOperator{\Fitt}{Fitt}

\DeclareMathOperator{\Jac}{Jac}
\DeclareMathOperator{\gr}{gr}

\DeclareMathOperator{\Sym}{Sym}
\DeclareMathOperator{\Tor}{Tor}

\newcommand{\myindent}{2.5em}
\makeatletter
\newtheoremstyle{myrmk}{0.5\topsep}{0.5\topsep}{
    \addtolength{\@totalleftmargin}{\myindent}
    \addtolength{\linewidth}{-\myindent}
    \parshape 1 \myindent \linewidth
}{-\myindent}{\bfseries}{.}{5pt plus 1pt minus 1pt}{}
\newtheoremstyle{myex}{\topsep}{\topsep}{
    \addtolength{\@totalleftmargin}{\myindent}
    \addtolength{\linewidth}{-\myindent}
    \parshape 1 \myindent \linewidth
}{-\myindent}{\bfseries}{.}{5pt plus 1pt minus 1pt}{}
\makeatother

\newtheorem{thm}{Theorem}[section]
\newtheorem{prop}[thm]{Proposition}
\newtheorem{lemm}[thm]{Lemma}
\newtheorem{coro}[thm]{Corollary}

\theoremstyle{definition}
\newtheorem{defn}[thm]{Definition}

\theoremstyle{myrmk}
\newtheorem{rmk}[thm]{Remark}

\theoremstyle{myex}
\newtheorem{examp}[thm]{Example}

\newcommand{\introemph}[1]{{\bf #1}}
\newcommand{\bodyemph}[1]{{\bf #1}}
\newcommand{\defemph}[1]{{\bf #1}}

\title{\textsc{Milnor Fiber Consistency via Flatness}}
\author{\sc Alex Hof}
\date{\begin{minipage}{\linewidth}
        \small{\bf Abstract.} We describe a new algebro-geometric perspective on the study of the Milnor fibration and, as a first step toward putting it into practice, prove powerful criteria for a deformation of a holomorphic function germ to admit a stratification on its domain partially satisfying the Thom condition and, more generally, to respect the Milnor fibration of the original germ in an appropriate sense. As corollaries, we obtain a method of partitioning the space of homogeneous polynomials of a fixed degree into finitely many locally closed subsets such that the fiber diffeomorphism type of the Milnor fibration is constant along each subset and a criterion under which deformations of a function with critical locus a complete intersection will be well-behaved.
\end{minipage}}

\begin{document}
\maketitle

\section{Introduction} \label{sec:intro}

Consider a nonzero holomorphic function germ $f: (\bC^{n+1}, 0) \to (\bC, 0)$ for an integer $n \ge 0$. Then the fiber over the origin in $\bC$ of $f$ is well-defined as a complex-analytic space germ $(V(f), 0)$ at the origin in $\bC^{n+1}$, which may in general be singular. There is also a well-defined local {\it smooth} fiber of $f$, given by the intersection of a sufficiently small ball around the origin in $\bC^{n+1}$ with the inverse image under $f$ of a sufficiently small nonzero value in $\bC$. Indeed, these nearby fibers fit together into a well-defined smooth locally trivial fibration called the \introemph{Milnor fibration}, which was introduced in \cite{milnor68} and has been much studied in the subsequent decades. For an introduction to the subject, see Chapter 3 of \cite{dimca92}, Chapter 9 of \cite{max_book}, or the handbook entry \cite{handbook_mf}; for more general classes of functions where analogues of the Milnor fibration exist, see, e.g., \cite{looijenga} or \cite{pt_irreg,art_fibr,art_sphere}.

There is much about the Milnor fibration which is not yet understood. In particular, given $f$, it is far from clear in general how to determine the structure of the Milnor fiber even at the level of homology, although there are techniques rendering the question tractable in various special cases. Perhaps the most compelling of these is one which has been known since the beginning --- Milnor showed in \cite{milnor68} that, if $(V(f), 0)$ has an isolated singularity at the origin, then the Milnor fiber of $f$ has the homotopy type of a bouquet of $n$-spheres, whose number, called the \introemph{Milnor number}, may be computed as the length of the quotient of the local ring by the Jacobian ideal $J_f$ of $f$.

Less powerfully but more generally, Kato and Matsumoto established in \cite{katmat73} that, if $s$ is the complex dimension of the singular locus of $(V(f), 0)$, then the Milnor fiber of $f$ is at least $(n-s-1)$-connected --- in particular, the reduced homology of the fiber is zero outside the degrees $[n-s, n]$ (for a perspective on this homological vanishing using the machinery of perverse sheaves see, e.g., \cite{max_van} or \cite{mpt_van}). Together with Milnor's stronger result in the case $s = 0$, this control of the Milnor fiber by the singularities of $(V(f), 0)$ is suggestive --- it now seems reasonable to wonder whether it is possible to compute the Milnor fiber entirely from the critical locus of $f$ at $0$, endowed with the complex-analytic scheme structure given by considering it as the vanishing of the Jacobian ideal $J_f$.

Note that, for there to be any hope of an affirmative answer to this question, the scheme-theoretic perspective on the critical locus is essential --- all isolated singularities, for example, have the same singular locus set-theoretically, but of course a variety of Milnor fibers are possible. The aim of this paper is to give a first step toward the study of the Milnor fibration from this viewpoint.

As is often the case in these sorts of situations, it is difficult to obtain direct results of the kind we want about a given function germ on its own, but studying how the objects of interest vary in families is more feasible. Therefore we will consider holomorphic function germs $F: (\bC^{n+1} \times \bC^u, 0 \times 0) \to (\bC, 0)$, regarding the coordinates of $\bC^u$ as parameters describing a germ of a family of holomorphic maps $\bC^{n+1} \to \bC$. For the sake of convenience, we denote by $f: (\bC^{n+1}, 0) \to (\bC, 0)$ the germ given by restricting $F$ to $(\bC^{n+1} \times 0, 0 \times 0)$.

We are interested in the following question, stated somewhat vaguely: When does the Milnor fibration vary consistently in the family $F$ and, in particular, how can we control this behavior algebraically? In the case where $f$ defines an isolated singularity, it is not difficult to show that in an appropriate sense the diffeomorphism type of the semi-local smooth fiber remains fixed in the family. This gives rise to a method for computing the Milnor fiber of $f$ known as \introemph{Morsification}, wherein $f$ is perturbed slightly to have only Morse singularities, the Milnor fiber of each of which contributes a separate cycle to the homology of the semi-local smooth fiber --- see, e.g., Chapter 2 of \cite{agzv_ii} for the details of this method, or see Chapter 7 of \cite{looijenga} for a treatment in the more general setting of isolated complete intersection singularities.

In the non-isolated case, the situation is much more complicated, and there have been various attempts to establish algebraic conditions under which some sort of consistency can be obtained. One example, of a somewhat different flavor than our approach here, is the work of Massey on \introemph{L\^{e} numbers} (for which, see, e.g., \cite{massey_lebook,handbook_le}), coordinate-dependent quantities defined for an arbitrary hypersurface singularity which in some sense generalize the Milnor number in the isolated case. The constancy of these numbers at the origin in a family of holomorphic functions guarantees the constancy of the Milnor fiber at least on the homological level, and in stronger senses if the critical locus has small enough dimension; however, the fact that this constancy applies to the Milnor fiber itself rather than the semi-local smooth fiber means the condition is rather stringent, precluding, for example, splittings into simpler singularities as in the Morsification method.

The other main algebraic method for finding perturbations which give useful information about the Milnor fibration is to recover the finite determinacy present in the isolated case by restricting to families within an ideal with respect to which $f$ has \introemph{finite extended codimension} --- geometrically, such deformations keep the positive-dimensional parts of the critical locus more or less fixed but allow us to manipulate zero-dimensional features (e.g., by splitting off Morse points). This strain of research was initiated by Siersma in \cite{siersma_linesing}, which was the first of a series of papers \cite{siersma_icis, siersma_planecurve, de_jong, pell_series, pell_1d, zah_2d, nem_2d, bm_3d} by various authors using these and related algebraic techniques to obtain Milnor fiber information about functions with critical loci satisfying various restrictions on dimension, structure, transversal type, and so forth. There has also been a great deal of work refining and generalizing the machinery of smooth local fibrations used in these sorts of splitting-type results from a more geometric perspective --- see, for instance, the work of Tib\u{a}r and his collaborators in \cite{st_deform,pt_irreg,art_fibr,art_sphere,jt_img}. (After the completion of the present work, particularly relevant results in this vein have been obtained in the preprint \cite{jst_tame} and article \cite{cjt_const}, the first of which also contains a more algebraic \introemph{Jacobian criterion} allowing control of the Milnor fibration through consideration of the integral closure of the Jacobian ideal.)

The theoretical underpinnings of the finite determinacy approach introduced by Siersma were developed in the thesis \cite{pell_thesis} of his student Pellikaan and the subsequent articles \cite{pell_findet, pell_series, pell_1d} based on its contents; much later, Bobadilla expanded on them in \cite{bobadilla_findet,bobadilla_fec}, proving a general theorem on Morsification in this context, and used his results to provide answers to some open questions about equisingularity-related topics in \cite{bobadilla_qs}. Connections to other perspectives were also pursued around this time by Gaffney in \cite{gaffney_inv,gaffney_pairs}.

Despite these results, it seems unlikely that an approach based on finite determinacy will be enough to solve the riddle of the Milnor fibration --- as mentioned, such methods allow for manipulations of the critical locus only at isolated points, which puts most of the non-isolated phenomena we would like to be able to address outside their purview. Fortunately, the perspective that the Milnor fiber should be controlled by the scheme-theoretic critical locus yields a much more general approach which does not suffer the same limitations.

Specifically, if we believe that the critical scheme of $f$ should encode information about the Milnor fiber, then we should try to get the Milnor fibration to vary consistently in our family $F$ by imposing some kind of control on the critical loci of the functions in the family. The most obvious condition to require is \introemph{flatness} (Definition \ref{def:flatness}), ubiquitous in algebraic and complex-analytic geometry as the correct notion of consistency for families of objects like schemes or sheaves --- for the critical locus to determine the Milnor fibration {\it intrinsically}, we should hope that its flatness over the parameter space is enough to force the consistency of the smooth fibers, and in nice situations this does indeed turn out to be true. More generally, it seems that the information we want involves not just the critical scheme but also its inclusion into the ambient space, and so it is necessary to ask for consistency in an {\it embedded} sense --- that is, we should require the flatness over the parameter space of all $k$th-order \bodyemph{infinitesimal neighborhoods} of the critical locus in the ambient space, or equivalently of the \bodyemph{normal cone} (Definition \ref{def:nc}) to the critical locus in the ambient space.

Such a requirement will in fact turn out to be enough to guarantee the consistency of smooth fibers, as desired. Our first result in this direction, which we will prove in Section \ref{sec:stratthm}, ties the flatness of the normal cone of the family's critical locus (as given in Definition \ref{def:critloc}) to the theory of stratifications:

\begin{thm} \label{thm:strat}
    Let $u \ge 0$ be an integer, $X'$ a complex-analytic manifold, $X \subseteq X' \times \bC^u$ an open subset, $\pi: X \to \bC^u$ the restriction to $X$ of the projection onto $\bC^u$, and $F: X \to \bC$ a holomorphic function which is not constant on any connected component of a fiber of $\pi$.
    
    Then there exists a complex-analytic Whitney stratification of $X$ such that $X \setminus \Sigma_{F \times \pi}$ is the ambient stratum, the non-flat locus of the normal cone $C_{\Sigma_{F \times \pi}}X$ to the critical locus $\Sigma_{F \times \pi}$ in $X$ over $\bC^u$ is a union of strata, and the Thom ($a_{F \times \pi}$) condition with respect to the ambient stratum is satisfied on any stratum not contained in this non-flat locus.
\end{thm}

This theorem, as mentioned, will ultimately allow us to use the flatness of the normal cone to the critical locus to control the behavior of smooth fibers in a family of holomorphic function germs by way of a smooth locally trivial fibration similar to those considered in \cite{art_fibr} (although note Remark \ref{rmk:setgerm} for a technical distinction). Also as mentioned, it is not always necessary to consider the full normal cone --- the flatness of the critical locus will be enough when there are, in some appropriate sense, vector fields witnessing the existence of the Milnor fibration of the function being deformed. We collect both of these results in the following theorem, which we will prove in Section \ref{sec:fibthm}:

\begin{thm} \label{thm:fibr}
    Let $n, u \ge 0$ be integers, $F: (\bC^{n+1+u}, 0) \to (\bC, 0)$ a holomorphic function germ, and $\pi: (\bC^{n+1+u}, 0) \to (\bC^u, 0)$ the projection onto the last $u$ coordinates. Suppose $f := F|_{\bC^{n+1} \times 0}$ is non-constant.
    
    Suppose that either of the following conditions holds:
    \begin{enumerate}[label=(\roman*)]
        \item There exists a sufficiently small open set $U \ni 0$ in $\bC^{n+1} \times 0$ such that $f: U \to \bC$ is well-defined, and, for each critical point $p \in (U \setminus 0) \cap \Sigma_f$ of $f$ other than the origin, the following holds: There exists a germ $V$ of a holomorphic vector field at $p$ such that the derivative $Vf$ of $f$ along $V$ vanishes identically as a complex-analytic function germ at $p$ and $V(p)$ is not tangent to the sphere of radius $|p|$ centered at the origin.
        
        Moreover, the intersection with $\pi^{-1}(0) = (\bC^{n+1} \times 0, 0)$ of the non-flat locus of $(\Sigma_{F \times \pi}, 0)$ over $(\bC^u, 0)$ is contained in the origin.
        
        \item The intersection with $\pi^{-1}(0) = (\bC^{n+1} \times 0, 0)$ of the non-flat locus of the normal cone of $(\Sigma_{F \times \pi}, 0)$ in $(\bC^{n+1+u}, 0)$ over $(\bC^u, 0)$ is contained in the origin.
    \end{enumerate}
    Then, if we fix a small enough $\varepsilon > 0$, we can choose $\delta > 0$ and $\gamma > 0$ small enough relative to $\varepsilon$ so that the following hold: The representatives $F: B_\varepsilon \times B_\gamma \to \bC$ and $\Sigma_{F \times \pi} \subset B_\varepsilon \times B_\gamma$ exist and, if we let $\Delta := (F \times \pi)(\Sigma_{F \times \pi}) \subset \bC \times B_\gamma$ be the discriminant, the restriction $$F \times \pi: (B_\varepsilon \times B_\gamma) \cap (F \times \pi)^{-1}((D_\delta \times B_\gamma) \setminus \bar \Delta) \to (D_\delta \times B_\gamma) \setminus \bar \Delta$$ defines a smooth locally trivial fibration, where $B_\varepsilon \subset \bC^{n+1}$, $D_\delta \subset \bC$, and $B_\gamma \subset \bC^u$ are the open balls of the specified radii.
\end{thm}

Thus, if the scheme-theoretic critical locus of the family defined by $F$ is consistent over the parameter space in an appropriate sense, at least away from the origin, then the semi-local smooth fibers in the family fit together into a smooth locally trivial fibration. By considering the restriction over the origin in the parameter space, we see that the fiber of this fibration is simply the Milnor fiber of $f$ if the radius is taken small enough. Therefore, within such a family, it is possible to study the smooth fiber by taking a slight perturbation, as in the isolated case --- this is what we mean by ``consistency''. As mentioned, this gives an important first step toward the larger project of computing information about the Milnor fibration of $f$ from the invariants of the Jacobian ideal of $f$ and lends credibility to the idea that such a project will be feasible.

The rest of this paper is organized as follows. Section \ref{sec:prelim} gives short introductions to a few concepts and results which will be important in the remainder of the paper, with references to more complete treatments; experienced readers may wish to skip over some or all of them, depending on their areas of expertise. Section \ref{sec:lem} lays the algebraic groundwork for the proofs of the main theorems via a series of lemmas, which are then applied to prove Theorem \ref{thm:strat} in Section \ref{sec:stratthm} and Theorem \ref{thm:fibr} in Section \ref{sec:fibthm}.

Section \ref{sec:homog} explores consequences of the main theorems for the study of homogeneous polynomials --- in particular, Corollary \ref{cor:homstrat} gives a partition of the space of homogeneous polynomials of a fixed degree into finitely many locally closed subsets along which the fiber diffeomorphism type of the Milnor fibration is constant. Section \ref{sec:ci} treats the special case, studied by various prior authors, of functions with critical locus a complete intersection, for which our theorems' hypotheses become particularly simple.

Section \ref{sec:comp} compares the theorems' consequences to the existing techniques described above and other past results. Finally, Section \ref{sec:ex} collects concrete examples not included elsewhere in the paper, including an instance of the theorems' use to compute homological data.

\begin{description}
    \item[Acknowledgments.] The author is indebted to his advisor, Lauren\c{t}iu Maxim, for his guidance and support throughout this project. He is also grateful to Javier Fern\'{a}ndez de Bobadilla, Dirk Siersma, and Terence Gaffney for interesting conversations. Much of the work in this paper was carried out during a Research Assistantship supported by funds from the University of Wisconsin-Madison's 2021 Fall Competition.
    
    The author also owes thanks to two anonymous referees for numerous useful suggestions and improvements to the exposition, to one of the two for proposing the comparison of Remark \ref{rmk:thomirreg}, and to the other for suggesting the deformation used in Example \ref{ex:novec}.
\end{description}

\section{Preliminaries and References}
\label{sec:prelim}

\subsection{Complex-Analytic Geometry and Cones}

Here we recall some relevant parts of the theory of \bodyemph{complex-analytic spaces}; for an introduction to these objects and their properties the reader is encouraged to refer to both \cite{gunningrossi} and \cite{fischer}. Essentially, they are the complex-analytic analogues to finite-type separated schemes over $\bC$, capturing holomorphic behavior instead of algebraic. In particular, germs of complex-analytic spaces correspond directly to quotients of convergent power series rings (see Section 0.21 of \cite{fischer}), so on a local level the study of these spaces is nearly just the study of a certain kind of affine scheme (although a certain amount of care must be taken as soon as any kind of non-locality is introduced). For this reason, and for the sake of terminological convenience, we will occasionally refer to non-reduced behavior and other algebraic properties as ``\bodyemph{scheme-theoretic}'' rather than ``complex-analytic-space-theoretic'', as would be strictly correct.

As mentioned in Section \ref{sec:intro}, the natural non-reduced structure on a function's critical locus will be of primary importance for us. Here we give the formal definition:

\begin{defn}[e.g., \cite{teissier_hunt}] \label{def:critloc}
    Let $\Phi: X \to Y$ be a flat map set-theoretically of pure dimension $d$ between finite-type schemes over a perfect field or complex-analytic spaces, or a germ of such a map. Then we define the \defemph{Jacobian ideal sheaf of $\Phi$} to be $\cJ_\Phi := \Fitt_d(\Omega_\Phi)$, the $d$th \defemph{Fitting ideal sheaf} (for which see, e.g., Section 1 of \cite{teissier_hunt} or Section 20.2 of \cite{eisenbud}) of the \defemph{sheaf of relative K\"{a}hler differentials} of $\Phi$ (for which see, e.g., Section 2 of \cite{teissier_hunt} or Section II.8 of \cite{hartshorne}). In the case of a map-germ, we write $J_\Phi$ in place of $\cJ_\Phi$ and omit ``sheaf'' from the terminology.
    
    We define the \defemph{critical locus of $\Phi$} to be the vanishing in $X$ of $\cJ_\Phi$ (respectively $J_\Phi$) and denote it by $\Sigma_\Phi$.
\end{defn}

In the case of a map of complex-analytic manifolds, the local computation can be made quite explicitly:

\begin{examp}
    Let $n, k \ge 0$ be integers and $\Phi: (\bC^{n+k}, 0) \to (\bC^k, 0)$ a holomorphic map-germ such that $\dim \Phi^{-1}(0) = n$. Then $J_\Phi$ is precisely the ideal of $k \times k$ minors of the \defemph{Jacobian matrix} $$\Jac \Phi := \begin{bmatrix}
        \frac{\partial \Phi_1}{\partial x_1} & \frac{\partial \Phi_1}{\partial x_2} & \cdots & \frac{\partial \Phi_1}{\partial x_{n+k}} \\
        \vdots & \vdots & \ddots & \vdots \\
        \frac{\partial \Phi_k}{\partial x_1} & \frac{\partial \Phi_k}{\partial x_2} & \cdots & \frac{\partial \Phi_k}{\partial x_{n+k}}
    \end{bmatrix},$$ where $\Phi_1, \ldots, \Phi_k$ are the component functions of $\Phi$ and $x_1, \ldots, x_{n+k}$ are the coordinates of $\bC^{n+k}$.
\end{examp}

Beyond this, we will often be particularly interested in \bodyemph{complex-analytic cones} (see \cite{am_cones}) over other complex-analytic spaces --- roughly, these are given by maps of complex-analytic spaces whose fibers are affine cones over (possibly weighted) projective schemes, although the fully rigorous definition expresses this requirement in terms of diagrams to properly account for non-reduced behavior. In the spirit of Serre's GAGA principle (\cite{gaga}), we then have the following result, which will make these spaces easier to manipulate in practice:

\begin{prop}[\cite{am_cones}] \label{prop:conecorresp}
    Let $(X, \cO_X)$ be a complex-analytic space. Then there is an anti-equivalence of categories, given by taking the analytic spectrum, between the category of finitely-presented nonnegatively-graded sheaves of $\cO_X$-algebras with degree-zero part $\cO_X$ and the category of complex-analytic cones over $X$.
\end{prop}

For the construction of the analytic spectrum in this case, see \cite{houzel_specan}.

There are two cones which will be of especial importance to us in what follows. The first of these captures information about the embedding of a closed subspace into a complex-analytic space:

\begin{defn} \label{def:nc}
    Let $(X, \cO_X)$ be a complex-analytic space and $Y$ the closed subspace cut out by an ideal sheaf $\cI$. Then the \defemph{normal cone $C_Y X$ to $Y$ in $X$} is the complex-analytic cone over $Y$ corresponding under Proposition \ref{prop:conecorresp} to the associated graded sheaf of algebras $$\gr_\cI \cO_X := \bigoplus_{i \ge 0} \cI^i/\cI^{i+1},$$ which is finitely presented by Proposition 1.17 of \cite{am_cones} and the Oka Coherence Theorem (see Section 0.12 of \cite{fischer}).
\end{defn}

Geometrically, this can be thought of roughly as the ``straightening-out of the ambient space in directions normal to the subspace''; in the case where both spaces are smooth, it is simply the normal bundle to the subspace in the usual sense of differential geometry. There is also a corresponding notion of a normal cone of a closed subscheme in algebraic geometry (see, e.g., \cite{fulton_it}), in which, instead of the analytic spectrum, we have the typical relative spectrum; these notions agree in the sense that, if we apply each to germ of an inclusion of complex-analytic spaces, we get respectively the analytic spectrum and the usual spectrum of the same finitely-presented algebra over the appropriate local ring.

In dealing in particular with complete intersections, we will make use of the following fact, which appears as Exercise 17.16 of \cite{eisenbud} and is originally due to Rees:

\begin{prop}[\cite{rees}] \label{prop:cinc}
    An inclusion $Y \hookrightarrow X$ of locally Noetherian schemes is a \defemph{regular embedding} --- that is, is locally cut out by a regular sequence --- if and only if $C_Y X$ is a vector bundle over $Y$. Moreover, a given choice of local generators for the ideal sheaf of $Y$ in $X$ will form a regular sequence exactly when the generators' initial forms give a local basis of sections for $C_Y X$.
\end{prop}

The second cone we will be interested in is in some senses similar, but captures instead the cotangent directions perpendicular to the smooth points of fibers of a map between complex manifolds:

\begin{defn}
    Let $\Phi: (X, \cO_X) \to (Y, \cO_Y)$ be a map between smooth complex-analytic spaces which is a submersion on an open dense subset $U$ of $X$. Then the \defemph{relative conormal space} of $\Phi$ is the closed subset of the cotangent bundle $T^*X$ given by $$T_\Phi^*X := \overline{\{(x, \eta) \in T^*U \mid \eta(\ker D\Phi|_x) = 0\}}.$$
\end{defn}

For our purposes, the following characterization will be more useful. It is straightforward to see that, under the anti-equivalence of Proposition \ref{prop:conecorresp}, the cotangent bundle of $X$ corresponds to the symmetric algebra on the sheaf $\Theta_X$ of vector fields on $X$, and the pullback along $\Phi$ of the cotangent bundle of $Y$ is likewise given by the symmetric algebra on $\Phi^*\Theta_Y$. Then we have:

\begin{prop} \label{prop:relconorm}
    $T_\Phi^*X$ corresponds under Proposition \ref{prop:conecorresp} to the sheaf of $\cO_X$-algebras which is the image of the map $\Sym(D\Phi): \Sym(\Theta_X) \to \Sym(\Phi^*\Theta_Y)$ given by the Jacobian of $\Phi$.
\end{prop}

\begin{description}
    \item[Proof] Around points of $X$ where $\Phi$ is a submersion, working in local coordinates makes this immediate. To prove it in general, we localize around each point of $T^*X$ --- which we regard as the analytic spectrum of $\Sym(\Theta_X)$ --- and then use the fact that the germ of the closed subcone of the cotangent bundle corresponding to the algebra described above is cut out by a prime ideal in the local ring and apply Theorem III.C.16 of \cite{gunningrossi} to show that it does in fact give the closure of the part over the submersive locus.
\end{description}

For both of these cones, the primary questions we will be interested in revolve around the following notion, which was introduced by Serre in \cite{gaga}:

\begin{defn} \label{def:flatness}
    A map $R \to S$ of rings is said to be \defemph{flat} if tensor by $S$ over $R$ preserves exact sequences of modules. A map of schemes or complex-analytic spaces is said to be \defemph{flat} precisely when all of the induced maps on local rings are flat.
\end{defn}

Expositions in the algebraic setting can be found in most textbooks on commutative algebra and algebraic geometry --- see, e.g., \cite{matsumura}, \cite{eisenbud}, or \cite{hartshorne}. The complex-analytic perspective on the subject is covered in Chapter 3 of \cite{fischer} and in \cite{hironaka}, particularly the Appendix; additionally, a criterion for the flatness of complex-analytic map germs in terms of lifting relations between defining equations can be found in Chapter 1 of \cite{stevens}. We will discuss this criterion in more detail in Lemma \ref{lem:diffeq}.

\begin{rmk}
    As mentioned in Section \ref{sec:intro}, flatness is generally accepted as the correct notion of consistency for algebro-geometric or complex-analytic objects parameterized over a base. It is common to speak of the geometric meaning of this consistency as being somewhat mysterious and unintuitive --- however, this need not be the case. As Hironaka points out in the Appendix to \cite{hironaka} (or see conditions $(4), (4')$ of Theorem 22.3 of \cite{matsumura}), in the Noetherian setting flatness at a point is equivalent to the triviality of the normal cone to the fiber over the tangent cone to the base, a straightforwardly geometric condition.
\end{rmk}

We conclude with the following remark about the interplay between cones and their corresponding sheaves of algebras, and between schemes over $\bC$ and their analytifications (for which see \cite{gaga}):

\begin{rmk} \label{rmk:alg}
    In what follows, we will very often be interested in verifying the flatness of complex-analytic spaces and cones over other complex-analytic spaces. To use these notions freely, we would like to be able to say the following:
    \begin{itemize}
        \item Given a map $(X, \cO_X) \to (Y, \cO_Y)$ of complex-analytic spaces and a complex-analytic cone over $X$, the cone's flatness over $Y$ as a complex-analytic space is everywhere equivalent to that of the corresponding sheaf of $\cO_X$-algebras over $\cO_Y$.
        
        \item Given a map $(X, \cO_X) \to (Y, \cO_Y)$ of separated finite-type schemes over $\bC$, the flatness of $X$ over $Y$ is everywhere equivalent to that of $X$'s analytification over $Y$'s analytification.
    \end{itemize}
    Both of these turn out to be true, essentially because of the faithful flatness of all maps of the form $$\bC\{u_1, \ldots, u_a\}[v_1, \ldots, v_b]_{(u_1, \ldots, u_a, v_1, \ldots, v_b)} \to \bC\{u_1, \ldots, u_a, v_1, \ldots, v_b\},$$ where $\bC\{u_1, \ldots, u_a\}$ and $\bC\{u_1, \ldots, u_a, v_1, \ldots, v_b\}$ are the appropriate convergent power series rings. Thus we can address questions of flatness in the algebraic or partially algebraic setting without needing to work explicitly with the analytifications of the objects involved.
\end{rmk}

\subsection{Stratifications and the Milnor Fibration}

In dealing with singular spaces from a set-theoretic perspective, we often wish to apply the machinery of differential geometry, either by working with a smooth ambient space or by decomposing the singular one into smooth pieces. This latter approach gives rise to the theory of \bodyemph{stratifications}, and, in particular, when equisingularity along the strata is required, \bodyemph{Whitney stratifications}. Brief introductions to these topics can be found in Chapter 1 of \cite{dimca92} and Chapter 1 of Part I of \cite{smt}; for a much more comprehensive overview, see \cite{handbook_strat} or \cite{larabida}. Here we recall a few facts and definitions from stratification theory which will be particularly relevant in what follows and use them to briefly discuss the Milnor fibration, introductory references for which have already been provided in Section \ref{sec:intro}.

Since it is important to the statement of Theorem \ref{thm:strat}, we recall explicitly the following relative version of the Whitney ($a$) condition:

\begin{defn}
    Let $X$ be a closed subset of a smooth manifold $M$ and $f$ the restriction to $X$ of a differentiable map from $M$ to another smooth manifold $N$. We say that a stratification of $X$ satisfies the \defemph{Thom ($a_f$) condition} on a stratum $S$ with respect to another stratum $S'$ such that $S \subset \overline{S'}$ if $f$ has constant rank on both strata and the following condition holds: If $s$ is a point of $S$ and $s_i$ is a sequence of points of $S'$ approaching $s$ such that the limiting tangent space $T_f := \lim_{i \to \infty} T_{s_i} (S' \cap f^{-1}(f(s_i)))$ to fibers of $f$ exists, then $T_s (S \cap f^{-1}(f(s))) \subseteq T_f$. We say that the stratification satisfies this condition if it does so on each stratum, and that a stratification of a map of real- or complex-analytic spaces does so if it does so locally with respect to the local embeddings into Euclidean space.
\end{defn}

Some authors prefer to write this as ``$A_f$'' rather than ``($a_f$)''. Unlike Whitney stratifications, stratifications satisfying the Thom condition are not guaranteed to exist. In \cite{hironaka}, particularly Section 5, Hironaka posits and explores a nebulous relationship between the flatness of $f$ and the existence of such a stratification, obtaining, among other things, the following result:

\begin{prop}[\cite{hironaka}, Section 5, Corollary 1] \label{prop:af}
    Let $f: X \to Y$ be a nowhere-constant complex-analytic map with $Y$ smooth of complex dimension 1 such that $f$ has isolated critical values. Then there exists a complex-analytic Whitney stratification of $f$ which satisfies the Thom ($a_f$) condition.
\end{prop}

Hironaka's proof of this fact, which uses Nash-style modifications, is quite involved --- in other works, such as \cite{dmod_strat}, \cite{hms}, \cite{leteissier}, \cite{gaffney_pairs}, and \cite{gafrang}, it is more common to establish Thom conditions using relative conormal spaces, which will be our approach for Theorem \ref{thm:strat}.

The utility of Thom stratifications for proving fibration theorems such as Theorem \ref{thm:fibr} lies in the following well-known result, also due to Thom:

\begin{lemm}[Thom's first isotopy lemma \cite{mather, mather_pub}] \label{lem:tfil}
    Let $X$ be a closed subset of a smooth manifold $M$ and $f: M \to N$ a map to another smooth manifold $N$ whose restriction to $X$ is proper. Then, if $X$ has a Whitney stratification such that the restriction of $f$ to each stratum is a smooth submersion, $f$ is a stratified-homeomorphically locally trivial fibration whose trivializations are moreover diffeomorphic along each stratum.
\end{lemm}

For a long time the canonical reference on this result was Mather's lecture notes \cite{mather}; these were eventually published with slight revisions as \cite{mather_pub}. Among the lemma's uses are the proof of a topological local triviality statement for Whitney stratifications --- see, e.g., Theorem 4.2.17 of \cite{handbook_strat}. We mention it explicitly mainly to make the following note about the exact nature of the trivializations involved:

\begin{rmk}
    It is common to omit the statement about stratum-wise diffeomorphism entirely, and state the lemma in terms of topological triviality alone. However, as Goresky and MacPherson note in their statement of the lemma in Part I, Chapter 1, Section 5 of \cite{smt}, the local trivialization homeomorphisms are indeed smooth along each stratum --- in fact, the proof in Mather's notes shows that the inverses are smooth as well, so they are even stratum-wise diffeomorphisms. This is because the trivializations and their inverses arise from the integration of controlled vector fields, which are by definition smooth along each stratum and hence yield flows which are likewise stratum-wise smooth. Thus the proof of Thom's first isotopy lemma gives a stronger result than is usually explicitly stated, even in \cite{mather, mather_pub}.
\end{rmk}

In particular, Proposition \ref{prop:af} and Lemma \ref{lem:tfil} can be applied as in \cite{le} to prove the existence of the following form of the Milnor fibration, sometimes called the \bodyemph{Milnor tube fibration} or \bodyemph{Milnor-L\^{e} fibration}:

\begin{prop}[\cite{milnor68, le}] \label{prop:mf}
    Let $f: (\bC^{n+1}, 0) \to (\bC, 0)$ be a germ of a holomorphic function. Then, if we fix a small enough $\varepsilon > 0$, we can choose $\delta > 0$ small enough relative to $\varepsilon$ so that the restriction $$f: B_\varepsilon \cap f^{-1}(D_\delta^*) \to D_\delta^*$$ defines a diffeomorphically locally trivial fibration, where $B_\varepsilon \subset \bC^{n+1}$ denotes the open ball of radius $\varepsilon$ and $D_\delta^* \subset \bC$ denotes the punctured open disk of radius $\delta$. Moreover, this is independent of the $\varepsilon$ and $\delta$ involved, provided they are indeed taken sufficiently small, and of the chosen local coordinates.
\end{prop}

L\^{e}'s proof, which we will generalize to prove Theorem \ref{thm:fibr}, is essentially an application of Proposition \ref{prop:af} to obtain the transversality of smooth fibers with the boundary sphere of $B_\varepsilon$ necessary for the use of Lemma \ref{lem:tfil}. Lemma \ref{lem:tfil} can likewise be used to give the independence of the resulting fibration from the chosen $\varepsilon, \delta$; for a proof that it is independent of the chosen coordinate system more broadly, see, e.g., Proposition 6.2.12 of \cite{handbook_mf}.

\begin{rmk} \label{rmk:radius}
    The method of proof described above gives us an effective criterion to confirm a chosen $\varepsilon$ is sufficiently small to be used to define the Milnor fibration. Specifically, it is enough that $S_\varepsilon$ and all smaller spheres centered at the origin are transverse to the strata of a Whitney stratification of $f$ satisfying the ($a_f$) condition; Whitney's condition ($b$), applied using lines through the origin, is enough to guarantee that this is indeed possible in this case. (More generalizably, one can use the \bodyemph{Curve Selection Lemma} to establish this transversality --- see Proposition 3.9 of \cite{hironaka_suban}.)
    
    It should be noted that this is a sufficient condition, but not a necessary one. In other circumstances, where it may not be possible to obtain Whitney stratifications of this type, the existence and independence of chosen radii of analogues to the Milnor fibration can also be established using more relaxed versions of these hypotheses --- see, e.g., Proposition 4.2 of \cite{art_fibr}.
\end{rmk}

\section{Algebraic Lemmas}
\label{sec:lem}

Here we state and prove various results, essentially algebraic in nature, which will be necessary for the proofs of Theorems \ref{thm:strat} and \ref{thm:fibr}.

\subsection{Normal Cones and Infinitesimal Neighborhoods}

To begin with, we discuss in more detail the properties of normal cones. We claimed in Section \ref{sec:intro} that the two notions of consistency of a closed embedding given by considering the infinitesimal neighborhoods and normal cone respectively were equivalent. We now prove this:

\begin{lemm} \label{lem:infnbhds}
    Let $\pi: X \to Y$ be a map of schemes (resp., complex-analytic spaces) and $Z \subseteq X$ a closed subscheme (resp., complex-analytic subspace). Denote by $\cI$ the sheaf of ideals corresponding to $Z$. Then the following are equivalent:
    \begin{itemize}
        \item For each $k \ge 0$, the $k$th-order infinitesimal neighborhood of $Z$ in $X$ (that is, the closed subspace of $X$ given by the vanishing of $\cI^{k+1}$) is flat over $Y$.
        
        \item The normal cone $C_Z X$ to $Z$ in $X$ is flat over $Y$.
    \end{itemize}
\end{lemm}

\begin{description}
    \item[Proof] As is noted in Remark \ref{rmk:alg} in the complex-analytic case and is immediate in the algebraic case, the flatness of the normal cone is equivalent to that of the associated graded sheaf of algebras $\gr_\cI \cO_X$. Since the tensor product distributes over direct sums, the flatness of the normal cone over $Y$ is hence equivalent to that of $\cI^k/\cI^{k+1}$ over $\cO_Y$ for all $k \ge 0$, which is to say of the flatness of $(\cI^k/\cI^{k+1})_p$ over $R := \cO_{Y,\pi(p)}$ for all $k \ge 0$ at each point $p$ of $X$. Hence we must show that this is equivalent to the flatness of $(\cO_X/\cI^{k+1})_p$ over $R$ for each $k \ge 0$ at each such point; hence we fix such a $p$ for the remainder of the proof.
    
    In general, the flatness of a module $M$ over $R$ is equivalent to the vanishing of $\Tor_1^R(M, N)$ for all $R$-modules $N$ and hence of $\Tor_i^R(M, N)$ for all integers $i \ge 1$ and all $R$-modules $N$. We can then prove each direction of the equivalence by applying the long exact sequence in Tor to the short exact sequences $$0 \to (\cI^k/\cI^{k+1})_p \to (\cO_X/\cI^{k+1})_p \to (\cO_X/\cI^k)_p \to 0$$ for $k \ge 0$; this gives us the exactness of the sequences $$\Tor_2^R\left(\left(\tfrac{\cO_X}{\cI^k}\right)_p, N\right) \to \Tor_1^R\left(\left(\tfrac{\cI^k}{\cI^{k+1}}\right)_p, N\right) \to \Tor_1^R\left(\left(\tfrac{\cO_X}{\cI^{k+1}}\right)_p, N\right) \to \Tor_1^R\left(\left(\tfrac{\cO_X}{\cI^k}\right)_p, N\right)$$ for any $R$-module $N$. It is then immediate that the flatness of the infinitesimal neighborhoods implies that of the normal cone, and the reverse implication can be proven by induction on $k$.
\end{description}

Hence we can pass freely back and forth between the two notions --- we will see that the flatness of the infinitesimal neighborhoods seems to arise more naturally as a desirable condition in proofs of other facts, whereas the flatness of the normal cone is much easier to establish by computations insofar as it involves a single finitely-presented algebra instead of a countably infinite family of modules.

In dealing with the normal cone, it is often useful to know that it will behave nicely under pullback. The following result uses the prior one to give sufficient conditions for this to occur:

\begin{lemm} \label{lem:ncpull}
    Let $X \to Y$ be a map of schemes (resp., complex-analytic spaces) and $Z \subseteq X$ a closed subscheme (resp., complex-analytic subspace). Let $\phi: Y' \to Y$ be a map of schemes (resp., complex-analytic spaces) and suppose that either of the following holds:
    \begin{enumerate}[label=(\roman*)]
        \item $\phi$ is a flat morphism.
        \item $C_Z X$ is flat over $Y$.
    \end{enumerate}
    Then $\phi^*(C_Z X) = C_{\phi^*Z} (\phi^*X)$.
\end{lemm}

\begin{description}
    \item[Proof] Let $\cI$ be the sheaf of ideals on $X$ corresponding to $Z$ and $\cI'$ the one on $X' := \phi^*X$ corresponding to $Z' := \phi^*Z$. By the anti-equivalence of Proposition \ref{prop:conecorresp}, it is enough to show that the natural map $\phi^*(\gr_\cI \cO_X) \to \gr_{\cI'} \cO_{X'}$ of sheaves of algebras on $X'$ is an isomorphism. This can be checked separately on each graded piece --- that is, it is enough to show that the natural map $\phi^*(\cI^k/\cI^{k+1}) \to (\cI')^k/(\cI')^{k+1}$ is an isomorphism.
    
    Now note that the formation of the infinitesimal neighborhoods of $Z$ in $X$ commutes with pullback --- that is, for each $k \ge 0$, $\phi^*(\cO_X/\cI^k) \cong \cO_{X'}/(\cI')^k$. Hence, by considering for each $k \ge 0$ the exact sequence $$0 \to \cI^k/\cI^{k+1} \to \cO_X/\cI^{k+1} \to \cO_X/\cI^k \to 0$$ and applying $\phi^*$, we see that it is enough to verify that the map $\cI^k/\cI^{k+1} \to \cO_X/\cI^{k+1}$ remains injective on pullback. This is immediate if our Condition (i) holds; if Condition (ii) holds, it follows from the long exact sequence in Tor by Lemma \ref{lem:infnbhds}.
\end{description}

\begin{rmk} \label{rmk:algc}
    In particular, Lemma \ref{lem:ncpull} allows us to compute normal cones {\it algebraically} when we are dealing with complex-analytic spaces defined by polynomials --- by applying the result together with the faithful flatness of the appropriate maps of local rings as in Remark \ref{rmk:alg}, we find that the formation of the normal cone will commute with analytification.
\end{rmk}

This lemma is, of course, also useful in the complex-analytic context, where it allows us to manipulate normal cones more freely in general.

As one example of the aforementioned tendency of the flatness of the infinitesimal neighborhoods to show up more readily in theoretical situations, we give the following special case of the flatness criterion of \cite{stevens}, which links the flatness of infinitesimal neighborhoods of a family's critical locus to the study of the partial differential equations satisfied by the functions in the family and plays a role in the proof of Theorem \ref{thm:fibr}:

\begin{lemm} \label{lem:diffeq}
Let $n, u \ge 0$ be integers, $F: (\bC^{n+1+u}, 0) \to (\bC, 0)$ a holomorphic function germ, and $\pi: (\bC^{n+1+u}, 0) \to (\bC^u, 0)$ the projection onto the last $u$ coordinates. Suppose $f := F|_{\bC^{n+1} \times 0}$ is non-constant.

Then, for each integer $d \ge 1$, the $(d-1)$st-order infinitesimal neighborhood of $(\Sigma_{F \times \pi}, 0)$ in $(\bC^{n+1+u}, 0)$ is flat over $(\bC^u, 0)$ if and only if every homogeneously first-order homogeneous degree-$d$ partial differential equation $$\sum_{\alpha_0 + \cdots + \alpha_n = d} c_{\alpha_0, \ldots, \alpha_n} \left(\frac{\partial f}{\partial x_0}\right)^{\alpha_0} \cdots \left(\frac{\partial f}{\partial x_n}\right)^{\alpha_n} = 0$$ satisfied by $f$ (with $c_{\alpha_0, \ldots, \alpha_n}: (\bC^{n+1}, 0) \to \bC$ holomorphic function germs) can be extended to a family $$\sum_{\alpha_0 + \cdots + \alpha_n = d} C_{\alpha_0, \ldots, \alpha_n} \left(\frac{\partial F}{\partial x_0}\right)^{\alpha_0} \cdots \left(\frac{\partial F}{\partial x_n}\right)^{\alpha_n} = 0$$ of such equations satisfied by $F$ (with $C_{\alpha_0, \ldots, \alpha_n}: (\bC^{n+1+u}, 0) \to \bC$ holomorphic function germs which restrict to the $c_{\alpha_0, \ldots, \alpha_n}$ on $(\bC^{n+1} \times 0, 0)$).
\end{lemm}

\begin{description}
    \item[Proof] This is immediate from the characterization of flatness in Chapter 1 of \cite{stevens}, which goes as follows. For a map $\phi: (X, 0) \to (S, 0)$ of germs of complex-analytic spaces and a chosen embedding $(\phi^{-1}(0), 0) \hookrightarrow (\bC^N, 0)$, we can write $\phi$ as the composition of an embedding $(X, 0) \hookrightarrow (\bC^N \times S, 0)$ with the projection to $(S, 0)$, as shown in Section 0.35 of \cite{fischer}. Then, if we let $G_1, \ldots, G_k$ be generators of the ideal of $\cO_{\bC^N \times S, 0}$ cutting out $(X, 0)$ in $(\bC^N \times S, 0)$ and $g_1, \ldots, g_k$ be their images under the quotient map to $\cO_{\bC^N \times 0, 0}$, $\phi$ is flat if and only if the following holds: Every relation $\sum_{i=1}^k r_ig_i = 0$ in $\cO_{\bC^N \times 0, 0}$ lifts to a relation $\sum_{i=1}^k R_iG_i = 0$ in $\cO_{\bC^N \times S, 0}$.
    
    In our case, $(X, 0)$ is the $(d-1)$st infinitesimal neighborhood of $(\Sigma_{F \times \pi}, 0)$ in $(\bC^{n+1+u}, 0)$, $(S, 0)$ is the parameter space $(\bC^u, 0)$, and $\phi$ is the projection $\pi$. Then, taking $N = n + 1$, we already have an embedding $(X, 0) \hookrightarrow (\bC^{n+1+u}, 0) \cong (\bC^N \times S, 0)$, and our $G_1, \ldots, G_k$ are precisely the degree-$d$ monomials in the partial derivatives of $F$. The result follows.
\end{description}

We will not require the full power of this lemma --- the $d=1$ case will be sufficient to demonstrate the facts about vector fields necessary for Condition (i) of Theorem \ref{thm:fibr} --- but the perspective it offers on the meaning of the flatness of the normal cone to the critical locus is intriguing and potentially worthy of further study in the future.

\subsection{Relative Conormal Spaces}

To prove Theorem \ref{thm:strat}, we will relate the flatness of the normal cone to properties of the relative conormal space, which can be used to track the limiting tangent planes relevant to the Thom condition. Recall by Proposition \ref{prop:relconorm} that the relative conormal of a map $\Phi: X \to Y$ of complex manifolds is given by the image of the map $\Sym(D\Phi): \Sym(\Theta_X) \to \Sym(\Phi^*\Theta_Y)$ of sheaves of algebras, and note in particular that $\Sym(\Theta_X)$ and $\Sym(\Phi^*\Theta_Y)$ are flat over $X$.

It then behooves us to establish the following general result, which says essentially that the behavior on pullback of a map of flat objects should be controlled by its cokernel:

\begin{lemm} \label{lem:spectral}
    Let $R$ be a ring, $0 \to A_3 \to A_2 \to A_1 \to A_0 \to 0$ an exact sequence of $R$-modules such that $A_1$ and $A_2$ are flat over $R$, and $B$ another $R$-module. Then the following hold:
    \begin{itemize}
        \item $H_2(A_\bullet \otimes_R B) \cong \Tor_1^R(A_0, B)$.
        \item $H_3(A_\bullet \otimes_R B) \cong \Tor_2^R(A_0, B)$.
        \item $\Tor_i^R(A_3, B) \cong \Tor_{i+2}(A_0, B)$ for all integers $i \ge 1$.
    \end{itemize}
\end{lemm}

\begin{description}
    \item[Proof] Let $\cdots \to F_2 \to F_1 \to F_0 \to B \to 0$ be a flat resolution of $B$. We can then form a first-quadrant double complex $A_\bullet \otimes_R F_\bullet$ whose $(i, j)$th entry is $A_i \otimes F_j$ using the usual sign trick (see, e.g., Section 2.7 of \cite{weibel}) and consider both of the corresponding spectral sequences (e.g., Section 5.6 of \cite{weibel}). By the flatness of the $F_j$ and the exactness of $A_\bullet$, the rows of this double complex are exact, so the convergence of the spectral sequence given by the filtration by rows tells us that the homology of the associated total complex vanishes.
    
    On the other hand, if we take the spectral sequence given by the filtration by columns, we find that the $E^1$ page has the following form:
    
    \begin{center}\begin{tikzcd}[column sep=small]
        & \vdots & \vdots & \vdots & \vdots & \\
        0 & \Tor_2(A_0, B) \ar[l] & 0 \ar[l] & 0 \ar[l] & \Tor_2(A_3, B) \ar[l] & 0 \ar[l] \\
        0 & \Tor_1(A_0, B) \ar[l] & 0 \ar[l] & 0 \ar[l] & \Tor_1(A_3, B) \ar[l] & 0 \ar[l] \\
        0 & \Tor_0(A_0, B) \ar[l] & \Tor_0(A_1, B) \ar[l] & \Tor_0(A_2, B) \ar[l] & \Tor_0(A_3, B) \ar[l] & 0 \ar[l] \\
    \end{tikzcd}\end{center}
    
    By considering the next few pages and using the convergence of this spectral sequence to the homology of the associated total complex, which we have already demonstrated is trivial, we obtain the desired isomorphisms.
\end{description}

There has been a fair amount of interest in the literature in results which characterize the behavior of the relative conormal space upon specialization, for example to the fiber over a point in the base --- see, e.g., Corollary 4.2.1 of \cite{hms}, Proposition 1.2.6 of \cite{leteissier}, Theorem 1.2.1 of \cite{fronces}, and Lemma 1.1 of \cite{gafmas}. Lemma \ref{lem:spectral} gives us control of this behavior on specialization along maps to the parameter space of a family using the flatness of the normal cone of the critical locus:

\begin{lemm} \label{lem:conorm_spec}
    Let $n, u \ge 0$ be integers, $F: (\bC^{n+1+u}, 0) \to (\bC, 0)$ a holomorphic function germ, and $\pi: (\bC^{n+1+u}, 0) \to (\bC^u, 0)$ the projection onto the last $u$ coordinates. Suppose $f := F|_{\bC^{n+1} \times 0}$ is non-constant.
    
    Suppose that the normal cone to $(\Sigma_{F \times \pi}, 0)$ in $(\bC^{n+1+u}, 0)$ is flat over $(\bC^u, 0)$. Then the restriction over $0 \in \bC^u$ of the relative conormal space of $F \times \pi$ satisfies $T_{F \times \pi}^*\bC^{n+1+u}|_{\pi^{-1}(0)} \cong T_f^*\bC^{n+1} \times \bC^u$ compatibly with the identification $T^*\bC^{n+1+u}|_{\pi^{-1}(0)} \cong T^*\bC^{n+1} \times \bC^u$.
\end{lemm}

\begin{description}
    \item[Proof] By using Proposition \ref{prop:relconorm} and considering the structure of the Jacobian matrix in this case, we find that the relative conormal space in this case is given as a cone over the space germ $(\bC^{n+1+u}, 0)$ by the graded algebra which is the image of the map $$\tilde\phi: \cO_{\bC^{n+1+u},0}[\xi_0, \ldots, \xi_n, \upsilon_1, \ldots, \upsilon_u] \to \cO_{\bC^{n+1+u},0}[\tau, \upsilon_1, \ldots, \upsilon_u]$$ which sends each $\xi_i$ to $\tfrac{\partial F}{\partial x_i} \tau$.
    
    Now, we can see that this map is in fact given by applying the functor $- \otimes_{\cO_{\bC^{n+1+u},0}} \cO_{\bC^{n+1+u},0}[\upsilon_1, \ldots, \upsilon_u]$, which corresponds under Proposition \ref{prop:conecorresp} to the product over $(\bC^{n+1+u}, 0)$ with a trivial rank-$u$ vector bundle, to the corresponding map $$\phi: \cO_{\bC^{n+1+u},0}[\xi_0, \ldots, \xi_n] \to \cO_{\bC^{n+1+u},0}[\tau].$$ Likewise, if we regard $\cO_{\bC^{n+1}, 0}$ as a quotient of $\cO_{\bC^{n+1+u},0}$ by the coordinate functions of $(\bC^u, 0)$, we find that $T_f^* \bC^{n+1}$ is given by the image of the induced map $\phi \otimes_{\cO_{\bC^{n+1+u},0}} \cO_{\bC^{n+1}, 0}$ in the quotient ring.
    
    It thus suffices to show that the formation of the image in this case commutes with passage to the quotient ring. Regarding $\phi$'s domain and codomain as graded $\cO_{\bC^{n+1+u},0}$-modules, we see that $$\coker \phi \cong \bigoplus_{i \ge 0} \frac{\cO_{\bC^{n+1+u},0}}{{J_{F \times \pi}}^i} \tau^i.$$ By Lemma \ref{lem:infnbhds}, we can then conclude that the flatness of the normal cone over $(\bC^u, 0)$ implies that of $\coker \phi$. This is to say that, for every $\cO_{\bC^u, 0}$-module $M$ and $i \ge 1$, $\Tor_i^{\cO_{\bC^u, 0}}(\coker \phi, M) = 0$.
    
    Observing that $\cO_{\bC^u, 0} \to \cO_{\bC^{n+1+u},0}$ is flat and hence that $\cO_{\bC^{n+1+u},0}[\xi_0, \ldots, \xi_n]$ and $\cO_{\bC^{n+1+u},0}[\tau]$ are flat $\cO_{\bC^u, 0}$-modules, we can then apply Lemma \ref{lem:spectral} to the exact sequence $$0 \to \ker \phi \to \cO_{\bC^{n+1+u},0}[\xi_0, \ldots, \xi_n] \xrightarrow{\phi} \cO_{\bC^{n+1+u},0}[\tau] \to \coker \phi \to 0$$ to find that it remains exact upon tensor over $\cO_{\bC^u, 0}$ by any module and hence that the image of the middle map in the sequence is preserved under such tensor operations. In particular, since tensor with the residue field of $\cO_{\bC^u, 0}$ gives us $\phi \otimes_{\cO_{\bC^{n+1+u},0}} \cO_{\bC^{n+1}, 0}$ in place of $\phi$, we have the claimed commutativity of image formation and thus the result follows.
\end{description}

This result will form the crux of our proof of Theorem \ref{thm:strat}.

\section{Proof of the Stratification Theorem}
\label{sec:stratthm}

We are now in a position to prove Theorem \ref{thm:strat}; once we have applied Lemma \ref{lem:conorm_spec}, our approach will follow the proof of the existence of Thom stratifications in \cite{dmod_strat} --- see also the proof of Proposition 8.3.10 of \cite{ks_som}, Proposition 10.6.5 of \cite{microhyp}, or Proposition 1.3.5 of \cite{leteissier}.

\begin{description}
    \item[Proof] We first observe that the non-flat locus of $C_{\Sigma_{F \times \pi}}X$ over $\bC^u$ in $X$ is indeed a complex-analytic subset of $X$. By, e.g., Theorem IV.9 of \cite{frisch}, the non-flat locus in $C_{\Sigma_{F \times \pi}}X$ itself is analytic. Moreover, since the complement of the zero section in $C_{\Sigma_{F \times \pi}}X$ is a $\bC^*$-bundle, and hence faithfully flat, over the projectivized normal cone $\bP C_{\Sigma_{F \times \pi}}X$, we can see that any failure of flatness at a point off the zero section of the normal cone will induce similar failures at every point in the same fiber of the bundle map and hence at the corresponding point of the zero section itself by the aforementioned analyticity. Thus, if we identify $\Sigma_{F \times \pi}$ with the zero section of the normal cone, we can see that the non-flat locus in $X$ we seek is simply its intersection with the non-flat locus in the normal cone itself, hence a complex-analytic subset of $X$ as claimed.
    
    Now consider the projectivization $\bP T_{F \times \pi}^* X$ of the relative conormal space. The projection $\bP T_{F \times \pi}^* X \to X$, being a proper complex-analytic map of complex-analytic spaces, admits a complex-analytic Whitney stratification by, e.g., Theorem 1 of Section 4 of \cite{hironaka}; the method of proof of that theorem shows that this can be taken such that $F \times \pi$ has constant rank on each stratum of $X$, the aforementioned non-flat locus is a union of strata, and $X \setminus \Sigma_{F \times \pi}$ is the ambient stratum --- this last condition can be met since $\bP T_{F \times \pi}^* X \to X$ is a $\bP^u$-bundle, and hence a surjective submersion, over $X \setminus \Sigma_{F \times \pi}$. We then claim that each stratum outside the non-flat locus satisfies the Thom ($a_{F \times \pi}$) condition with respect to this ambient stratum.
    
    To show this, consider such a stratum $S$ and observe that, since $S$ is contained in the critical locus of $F \times \pi$ and $\pi$ itself has full rank everywhere on $X$, the partial derivatives of $F$ with respect to the directions tangent to the fibers of $\pi$ must all vanish identically on $S$. Hence the fibers of $(F \times \pi)|_S$ are exactly those of $\pi|_S$ --- as such, we consider a point $q \in \bC^u$ and let $M = S \cap \pi^{-1}(q)$. This is a manifold because of our stipulation that $F \times \pi$ have constant rank on $S$; it remains to show that $M$'s tangent space at each point is contained in the intersection of all limiting tangent spaces to smooth points of fibers of $F \times \pi$. This is equivalent to the requirement that the conormal bundle $T_M^*X$ of $M$ in $X$ set-theoretically contain the restriction of the relative conormal space $T_{F \times \pi}^* X$ over $M$.
    
    Since $M$ is contained in $\pi^{-1}(q)$, we have a short exact sequence of vector bundles $$0 \to T_{\pi^{-1}(q)}^* X|_M \to T_M^* X \to T_M^* \pi^{-1}(q) \to 0,$$ which has a splitting given by the natural product decomposition of the cotangent bundle of $X' \times \bC^u$. Hence, since $T_{\pi^{-1}(q)}^* X|_M$ is a trivial vector bundle of rank $u$, we find that $T_M^* X \cong T_M^* \pi^{-1}(q) \times \bC^u$ compatibly with the identification $T^*X|_{\pi^{-1}(q)} \cong T^*\pi^{-1}(q) \times \bC^u$. Applying Lemma \ref{lem:conorm_spec} at each point, which is possible since we are outside the non-flat locus, then tells us that it is enough to show $T_f^*\pi^{-1}(q)|_M \subseteq T_M^* \pi^{-1}(q)$ set-theoretically, where $f := F|_{\pi^{-1}(q)}$.
    
    Let $T'$ be any stratum of $\bP T_{F \times \pi}^* X$ mapped (surjectively and submersively) to $S$. Denote by $X \times 0$ the zero section of the cotangent bundle, and define $T$ to be the inverse image of $T'$ under the $\bC^*$-bundle map $(T_{F \times \pi}^* X) \setminus (X \times 0) \to \bP T_{F \times \pi}^* X$. Since the projection $\bP T_{F \times \pi}^* X \to X$ is a proper stratified submersion, Thom's first isotopy lemma (\ref{lem:tfil}) tells us that it is smoothly locally trivial along each stratum. In particular, we can see that the inverse image $N$ of $M \subseteq S$ in $T$ is a complex-analytic manifold mapped surjectively and submersively to $M$.
    
    Now, since $M$ is contained in $\Sigma_{F \times \pi} \cap \pi^{-1}(q)$, we can see in particular that it is contained in a single fiber of $f$. By Theorem 3.4.2 of \cite{dmod_strat} (see also Theorem 3.1 of \cite{gafmas}), the restriction of $T_f^*\pi^{-1}(q)$ over this fiber of $f$ is set-theoretically the characteristic variety of a regular holonomic $D$-module and hence Lagrangian. In particular, it is isotropic. Thus, if we work around a point of $M$ and choose local coordinates $x_0, \ldots, x_n$ for $\pi^{-1}(q)$ such that $M$ is locally cut out by the vanishing of the coordinates $x_0, \ldots, x_{c-1}$, the vanishing of the canonical 1-form $\xi_0dx_0 + \ldots + \xi_ndx_n$ on $T_f^*\pi^{-1}(q)|_M$ gives us the vanishing of the 1-form $\xi_cdx_c + \ldots + \xi_ndx_n$ on $N \subset T_f^*\pi^{-1}(q)|_M$. Since $dx_c, \ldots, dx_n$ are linearly independent on $M$ and hence on $N$, this implies that $\xi_c, \ldots, \xi_n$ vanish on $N$, which is to say that $N \subseteq T_M^* \pi^{-1}(q)$.
    
    Repeating this reasoning for each such $T'$ --- and hence each such $N$ --- demonstrates that the complement of the zero section in $T_f^*\pi^{-1}(q)|_M$ is set-theoretically contained in $T_M^*\pi^{-1}(q)$. Since this complement is set-theoretically dense in $T_f^*\pi^{-1}(q)|_M$ and $T_M^*\pi^{-1}(q)$ is closed in $T^*\pi^{-1}(q)|_M$, the result follows.
\end{description}

We note some immediate comparisons between these objects and others considered elsewhere:

\begin{rmk} \label{rmk:nonmapstrat}
    What we obtain in Theorem \ref{thm:strat} is a stratification specifically of the space $X$, not of the map $F \times \pi$. That is, we do not produce a corresponding stratification of the codomain $\bC \times \bC^u$ such that each stratum of $X$ is mapped surjectively and submersively to one of its strata --- indeed, such a thing is not possible unless the critical locus of $F \times \pi$ is empty, since our stratification has its complement as the ambient stratum and $(F \times \pi)(X \setminus \Sigma_{F \times \pi}) = (F \times \pi)(X)$. In particular, although it is very common in the literature to consider stratifications of maps, often specifically those making a chosen fiber a union of strata, what we produce is not of this form.
\end{rmk}

\begin{rmk} \label{rmk:thomirreg}
    It is of interest to compare the non-flat locus of Theorem \ref{thm:strat} to other loci from the literature measuring the failure of the Thom condition in various senses. In particular, there is the \bodyemph{Thom irregularity locus} defined in \cite{pt_irreg}; this is not directly comparable, since it concerns the Thom condition with respect to all pairs of strata, not just the ambient one, but an analogous \bodyemph{partial Thom irregularity locus} which deals only with the ambient stratum appears in \cite{art_fibr} under the name ``$\partial$-Thom regularity set''. If we take a refinement of our stratification by intersecting with a particular fiber of $F \times \pi$, Theorem \ref{thm:strat} is then nearly enough to guarantee the containment of this locus in our non-flat locus; however, since the definition of \cite{art_fibr} involves a requirement for points to lie on {\it positive-dimensional} strata specifically, and our methods do not preclude the existence of a stratum outside the non-flat locus which has zero-dimensional intersection with some fiber, our arguments are not enough to establish this containment in general.
\end{rmk}

\section{Proof of the Fibration Theorem}
\label{sec:fibthm}

We can now prove Theorem \ref{thm:fibr} by applying either Lemma \ref{lem:diffeq} or Theorem \ref{thm:strat} to obtain the transversality necessary for Lemma \ref{lem:tfil}:

\begin{description}
    \item[Proof] Let $\Gamma, E > 0$ be small enough that $F$ is defined on $B_E \times B_\Gamma \subset \bC^{n+1} \times \bC^u$ and note that $\Sigma_{F \times \pi}$ thus has a well-defined representative as a closed complex-analytic subspace of $B_E \times B_\Gamma$. For any $0 < r < E$, let $S_r$ denote the sphere of radius $r$ centered at the origin in $\bC^{n+1}$. We will begin by arguing that, for small enough $0 < \varepsilon < E$, we can produce for each $p \in \Sigma_{F \times \pi} \cap (S_\varepsilon \times 0)$ a neighborhood $U_p$ of $p$ in $B_E \times B_\Gamma$ such that, for each $q \in U_p \setminus \Sigma_{F \times \pi}$, the fiber $(F \times \pi)^{-1}((F \times \pi)(q))$ meets $S_{|q|} \times B_\Gamma$ transversally at $q$.
    
    Suppose first that Condition (i) holds. Then, if we make $\varepsilon$ small enough, we have for each such $p$ a holomorphic vector field germ $V_p$ on $\bC^{n+1} \times 0$ at $p$ satisfying the stated hypotheses. Moreover, the condition tells us that $0 \times 0$ is either not contained in the non-flat locus of $\Sigma_{F \times \pi}$ over $B_\Gamma$ or is an isolated point of its intersection with $B_E \times 0$, so we can suppose by further shrinking $\varepsilon$ that $\Sigma_{F \times \pi}$ is flat over $B_\Gamma$ at each $p$. We can then apply Lemma \ref{lem:diffeq} in the case $d=1$ to the differential equation $V_pf = 0$ to obtain a germ at $p$ of a holomorphic vector field $\tilde V_p$ on $B_E \times B_\Gamma$ such that $\tilde V_pF = 0$, $\tilde V_p$ is everywhere tangent to the fibers of the projection $\pi$, and $\tilde V_p|_{B_E \times 0} = V_p$.
    
    By hypothesis, the tangent vector given by $V_p$ at $p$ is not contained in the tangent space $T_pS_\varepsilon$. By the openness of non-containment, we can then let $U_p$ be a neighborhood of $p$ in $B_E \times B_\Gamma$ small enough that $\tilde V_p$ has a representative defined on $U_p$ and the tangent vector given by this representative at each point $q$ of $U_p$ is not contained in the corresponding tangent space to $S_{|q|} \times B_\Gamma$. Since $\tilde V_p$ is everywhere tangent to the smooth fibers of $F \times \pi$, this $U_p$ thus has the desired transversality property.
    
    Now suppose instead that Condition (ii) holds. As before, the flatness hypothesis tells us that $0 \times 0$ is either not contained in the non-flat locus of the normal cone $C_{\Sigma_{F \times \pi}}(B_E \times B_\Gamma)$ over $B_\Gamma$ or is an isolated point of its intersection with $B_E \times 0$. We can then apply Theorem \ref{thm:strat} to obtain a Whitney stratification of $B_E \times B_\Gamma$ such that the non-flat locus is a union of strata and the Thom ($a_{F \times \pi}$) condition with respect to the ambient stratum $(B_E \times B_\Gamma) \setminus \Sigma_{F \times \pi}$ holds for the strata outside of it.
    
    Applying Whitney's condition ($b$) to lines through the origin (or using the Curve Selection Lemma) as in Remark \ref{rmk:radius}, we can then see that, for sufficiently small $0 < \varepsilon < E$, each point $p$ of $\Sigma_{F \times \pi} \cap (S_\varepsilon \times 0)$ lies on a stratum outside the non-flat locus which is transverse to $S_\varepsilon \times B_\Gamma$ at $p$. Using the fact that this stratum satisfies the Thom condition with respect to the ambient stratum, we find by the openness of non-containment that the smooth points $q$ of fibers of $F \times \pi$ in a sufficiently small neighborhood $U_p$ of $p$ are transverse to the corresponding submanifolds $S_{|q|} \times B_\Gamma$, as desired.
    
    Hence our hypotheses guarantee in all cases the existence of neighborhoods $U_p$ with the properties claimed. Now suppose toward a contradiction that we have a sequence of points of $S_\varepsilon \times B_\Gamma$ whose images under $F \times \pi$ are contained in the complement of the discriminant closure $\bar \Delta$ and approach $0 \times 0 \in \bC \times B_\Gamma$ such that, at each point of the sequence, the fiber of $F \times \pi$ fails to be transverse to $S_\varepsilon \times B_\Gamma$. By the compactness of $S_\varepsilon \times \bar B_\Gamma$, this sequence has a convergent subsequence, and we can see that the point $p$ to which it converges must lie in $S_\varepsilon \times 0$. Indeed, if $\varepsilon$ is small enough, we must have $p \in \Sigma_{F \times \pi}$ as well, since the spheres $S_\varepsilon$ will be eventually transverse to the regular locus of $(F \times \pi)^{-1}(0 \times 0)$ by, e.g., the Generic Whitney Lemma applied with Whitney's condition ($b$) as before and non-containment, as previously noted, is an open condition. Thus the definition of convergence implies that the sequence must be eventually contained in the corresponding neighborhood $U_p$. This gives the desired contradiction by the construction of $U_p$.
    
    Hence we can pick $\delta > 0$ and $0 < \gamma < \Gamma$ such that no such points have image in $D_\delta \times B_\gamma$. Then the restriction $$F \times \pi: (\bar B_\varepsilon \times B_\gamma) \cap (F \times \pi)^{-1}((D_\delta \times B_\gamma) \setminus \bar \Delta) \to (D_\delta \times B_\gamma) \setminus \bar \Delta$$ is a stratified submersion when we give the domain the stratification with strata $(B_\varepsilon \times B_\gamma) \cap (F \times \pi)^{-1}((D_\delta \times B_\gamma) \setminus \bar \Delta)$ and $(S_\varepsilon \times B_\gamma) \cap (F \times \pi)^{-1}((D_\delta \times B_\gamma) \setminus \bar \Delta)$. Moreover, it is proper, since it is a restriction over a subset of the target of the proper map $F \times \pi: \bar B_\varepsilon \times \bar B_\gamma \to \bC \times \bar B_\gamma$. As such, we can apply Thom's first isotopy lemma (\ref{lem:tfil}) to show that it is a locally trivial fibration both topologically and diffeomorphically along each stratum. Discarding the boundary stratum gives us the original result.
\end{description}

At this point it seems appropriate to address a technical subtlety in the statement of Theorem \ref{thm:fibr}:

\begin{rmk} \label{rmk:setgerm}
    Although we have avoided indicating the dependence notationally, the discriminant $\Delta$ in Theorem \ref{thm:fibr} is defined with reference to our choice of representatives and is not a priori independent of $\gamma, \delta, \varepsilon$ --- that is, it may not exist as a well-defined set germ. In particular, the statement of the theorem does not claim that the fibration obtained is independent of the chosen radii in general --- since our interest in this fibration is to study the Milnor fibration of the function being deformed rather than to come up with any kind of invariant of the deformation itself, such a result is not necessary for our purposes. For results on the existence of the discriminant as a well-defined complex-analytic set germ and the relationship between this property and the existence of fibrations independent of chosen radii, see \cite{jt_img} and \cite{art_fibr} --- in particular, note by Theorem 3.2 of \cite{jt_img} that $\Delta$ is a well-defined complex-analytic set germ in the case of a one-parameter deformation ($u = 1$), so we can always avoid such issues by restricting over a curve in our parameter space.
\end{rmk}

\section{Global Consequences}
\label{sec:homog}

If $f: \bC^{n+1} \to \bC$ is a homogeneous polynomial, the Milnor fibration has a global interpretation as well (see, e.g., Chapter 10 of \cite{max_book} or Example 3.1.9ff. of \cite{dimca92}). Specifically, the restriction $$f: \bC^{n+1} \setminus V(f) \to \bC^*$$ defines a smooth locally trivial fibration, called the \bodyemph{global} or \bodyemph{affine Milnor fibration} of $f$. It is not difficult to show that this is fiber diffeomorphic to the usual local Milnor fibration of $f$ at the origin.

Hence any homogeneous polynomial of degree $d$ has a well-defined smooth fiber, which comes with a monodromy action given by the action of a primitive $d$th root of unity on the polynomial's level sets, and these are together equivalent to the corresponding Milnor fibration information at the origin. We can then use Theorem \ref{thm:fibr} to obtain results about these global objects:

\begin{coro} \label{cor:homog}
    Let $n \ge 0$ and $d \ge 1$ be integers, $Y$ a connected complex-analytic space, and $$F: \bC^{n+1} \times Y \to \bC$$ a holomorphic map such that, for each point $y \in Y$, the map $f_y := F|_{\bC^{n+1} \times y}$ is a nonzero homogeneous degree-$d$ polynomial. Let $\pi: \bC^{n+1} \times Y \to Y$ be the projection and $\bP\Sigma_{F \times \pi}$ the projectivization of the critical locus of $F \times \pi$.
    
    Let $C_{\bP\Sigma_{F \times \pi}}(\bP^n \times Y)$ be the normal cone and suppose the natural projection $C_{\bP\Sigma_{F \times \pi}}(\bP^n \times Y) \to Y$ is flat. Then the fibrations over $\bC^*$ induced by the $f_y$ as $y$ varies throughout $Y$ are all fiber diffeomorphic to one another.
\end{coro}

\begin{description}
    \item[Proof] We first observe that, by pulling back along a resolution of singularities (see \cite{resoln}) $\tilde Y \to Y$, we may assume $Y$ is smooth --- this may result in a parameter space $\tilde Y$ with multiple connected components, but, since the fibration over any point of $\tilde Y$ is the same as the one over the corresponding point of the connected space $Y$, proving the constancy along each connected component will be sufficient to get consistency everywhere. Since the formation of the Jacobian ideal sheaf $\cJ_{F \times \pi}$ commutes with pullback along maps into $Y$ in general and the flatness hypothesis guarantees that the formation of the normal cone will do so as well by Lemma \ref{lem:ncpull}, the passage from $Y$ to $\tilde Y$ does not affect our hypotheses.
    
    Note that $\Sigma_{F \times \pi} \setminus (0 \times Y)$ is the pullback of $\bP\Sigma_{F \times \pi}$ along the natural $\bC^*$-bundle map $(\bC^{n+1} \times Y) \setminus (0 \times Y) \to \bP^n \times Y$ given by sending each point of $\bC^{n+1} \setminus 0$ to the point of $\bP^n$ corresponding to the line through the origin it lies on. In particular, since the bundle map is flat, the formation of the normal cone commutes with this pullback by Lemma \ref{lem:ncpull}, and hence the flatness of $C_{\bP\Sigma_{F \times \pi}}(\bP^n \times Y)$ over $Y$ implies that of $C_{\Sigma_{F \times \pi} \setminus (Y \times 0)}((\bC^{n+1} \times Y) \setminus (0 \times Y))$. Since flatness is a local property and the formation of the normal cone is a local process, this is exactly to say that $C_{\Sigma_{F \times \pi}}(\bC^{n+1} \times Y)$ is flat over $Y$ everywhere on $\bC^{n+1} \times Y$ except possibly at $0 \times Y$.
    
    Then, since the discriminant of $F \times \pi$ in this case will simply be $0 \times Y \subset \bC \times Y$, we can apply Theorem \ref{thm:fibr} to obtain, for each point $y \in Y$, a small neighborhood $N_y$ of $y$ and small enough values $\varepsilon_y, \delta_y > 0$ such that the maps $$f_{y'}: B_{\varepsilon_y} \cap {f_{y'}}^{-1}(D_{\delta_y}^*) \to D_{\delta_y}^*$$ as $y'$ varies in $N_y$ fit together into a smooth locally trivial fibration over $D_{\delta_y}^* \times N_y$. In particular, these maps define smooth locally trivial fibrations over the punctured disk which are all fiber diffeomorphic.
    
    However, the equivariance of homogeneous polynomials under scaling implies per Remark \ref{rmk:radius} that for a homogeneous polynomial every radius is a Milnor radius --- that is, every choice of $\varepsilon > 0$, no matter how large, can be used to define the Milnor fibration. In particular, the fibrations we have defined for the $f_{y'}$ in the ball $B_{\varepsilon_y}$ are, in fact, their Milnor fibrations at the origin, and hence fiber diffeomorphic to their affine Milnor fibrations. As such, the fiber diffeomorphism type of the affine Milnor fibration of $f_y$ is locally constant on $Y$, hence constant over all of $Y$ if we assume connectedness. As discussed, the result follows.
\end{description}

\begin{rmk}
    In this case, the Jacobian ideal sheaf $\cJ_{F \times \pi}$ is generated by the partial derivatives $\tfrac{\partial F}{\partial x_0}, \ldots, \tfrac{\partial F}{\partial x_n}$ and hence $\bP\Sigma_{F \times \pi}$ is given by the projective vanishing of these functions in $\bP^n \times Y$.
\end{rmk}

\begin{rmk}
    Note that here the homogeneity is crucial in guaranteeing that the smooth fibers we get from Theorem \ref{thm:fibr} as we vary the parameters are actually the Milnor fibers of the functions in the family. In general, this need not be true --- indeed, the fact that it is not is essential if we are to use the theorem to obtain splittings into simpler singularities as in the case of a Morsification of an isolated singularity or of our later Example \ref{ex:wup}.
\end{rmk}

We can apply this result to obtain a clearer picture of how the Milnor fibration varies over the space of all homogeneous polynomials:

\begin{coro} \label{cor:homstrat}
    Fix integers $n \ge 0$ and $d \ge 1$ and let $H_{n,d} \cong \bP^{\binom{n+d}{n}-1}$ be the space of degree-$d$ hypersurfaces in $\bP^n$, so that the $\binom{n+d}{n}$ projective coordinates of each point give the coefficients up to scaling of the monomial terms in a homogeneous polynomial defining the corresponding hypersurface. Then iteratively applying Corollary \ref{cor:homog} gives us a partition of $H_{n,d}$ into finitely many disjoint Zariski-locally-closed subsets so that the fiber diffeomorphism type of the affine Milnor fibrations of the corresponding defining polynomials is constant along each subset.
\end{coro}

\begin{description}
    \item[Proof] Let $\Sigma_{n,d}$ be the closed subscheme of $\bP^n \times H_{n,d}$ such that the fiber of $\Sigma_{n,d}$ over each point of $H_{n,d}$ is the scheme-theoretic singular locus of the corresponding projective hypersurface, and denote by $U_{n,d} \cong \bC^{\binom{n+d}{n}} \setminus 0$ the space of nonzero degree-$d$ homogeneous polynomials in $n+1$ variables. We will apply Corollary \ref{cor:homog} with the universal polynomial $$F: \bC^{n+1} \times U_{n,d} \to \bC$$ given by $F(x_0, \ldots, x_n, h) = h(x_0, \ldots, x_n)$; as usual, we take $\pi: \bC^{n+1} \times U_{n,d} \to U_{n,d}$ to be the projection. The projection $\bP\Sigma_{F \times \pi} \to U_{n,d}$ is then the pullback along the natural $\bC^*$-bundle map $U_{n,d} \to H_{n,d}$ of the projection $\Sigma_{n,d} \to H_{n,d}$ by dint of the fact that a homogeneous polynomial is contained in its own Jacobian ideal and so $\Sigma_{n,d}$ is defined by partial derivatives alone. Hence, since bundle maps are flat, Lemma \ref{lem:ncpull} tells us for any map $\psi: S \to H_{n,d}$ that pulling back the normal cone to $\psi^*\Sigma_{n,d}$ in $\bP^n \times S$ along the bundle map $\psi^*U_{n,d} \to S$ yields the normal cone to $\psi^*\bP\Sigma_{F \times \pi}$ in $\bP^n \times \psi^*U_{n,d}$. Thus, by Corollary \ref{cor:homog}, a partition of $H_{n,d}$ into finitely many disjoint Zariski-locally-closed subsets $S$ such that the normal cone to $\Sigma_{n,d}|_S$ in $\bP^n \times S$ is flat over $S$ will satisfy the desired constancy of the Milnor fibrations.
    
    To construct such a partition, we proceed as follows. Suppose $S$ is a closed subset of $H_{n,d}$. Then, by an appropriate modification (to account for the algebraicity of all objects involved) of the argument used in the proof of Theorem \ref{thm:strat}, the non-flat locus of the normal cone to $\Sigma_{n,d}|_S$ in $\bP^n \times S$ over $S$ is a closed algebraic subset of $\bP^n \times S$. Since the projection $\bP^n \times S \to S$ is projective and hence proper, the image of this locus in $S$ is also Zariski closed. Moreover, by generic flatness (e.g., Theorem 5.12 of \cite{FGAe} or Theorem 24.1 of \cite{matsumura}) applied to any irreducible component of $S$, this image cannot be all of $S$.
    
    As such, starting with $S = H_{n,d}$, we obtain the desired partition iteratively by letting $S'$ be this image of the non-flat locus, adding $S \setminus S'$ as one of the subsets in our partition, and then replacing $S$ with $S'$ and repeating the procedure. Since $H_{n,d}$ is Noetherian and the successive closed subsets $S$ form a descending chain, this procedure must terminate in finitely many steps. The result follows.
\end{description}

\begin{examp}
    For any $n \ge 0$, the space $\Sigma_{n, 2}$ of Corollary \ref{cor:homstrat} with $d=2$ will be defined by equations which are homogeneous linear in the variables $x_0, \ldots, x_n$. In particular, the partition of $H_{n,2}$ we obtain from our algorithm is simply the partition by the dimensions of the singular loci of the corresponding hypersurfaces, since $\Sigma_{n, 2}$'s restriction over each constant-dimensional locus is a projectivized vector bundle and so we can see by Proposition \ref{prop:cinc} that its normal cone will be a vector bundle over a projectivized vector bundle over the base, hence flat.
\end{examp}

In principle, Corollary \ref{cor:homstrat} can likewise be applied to obtain explicit partitions of $H_{n,d}$ for higher degrees $d$ --- however, in practice the number of variables involved makes direct computations expensive even for low $n \ge 2$ and $d \ge 3$.

\section{Complete Intersections}
\label{sec:ci}

One of the main restrictions which appears frequently in the works of Sierma's school is the hypothesis that the function to be deformed have a critical locus which is set-theoretically a complete intersection (\cite{siersma_linesing,siersma_icis,siersma_planecurve,pell_thesis,pell_findet,pell_series,pell_1d,zah_2d,nem_2d,gaffney_pairs}). Here we prove general results about the well-behavedness of deformations of such functions which are in an appropriate sense compatible with the complete intersection structure; compare these in particular with the results of Section 3 of \cite{zah_2d} and Section 4 of \cite{gaffney_pairs}.

To apply our techniques, of course, it is necessary to require some additional restriction on the scheme structure of the critical locus as well; we will discuss the case where the critical locus is {\it scheme-theoretically} a complete intersection, except possibly at the origin. Then the normal cone to the critical locus will behave well under deformations which respect the complete intersection structure:

\begin{prop} \label{prop:ci}
    Let $I = (g_1, \ldots, g_c) \subset \cO_{\bC^{n+1}, 0}$ be an ideal defining a germ at the origin of a complete intersection of codimension $c$ in $\bC^{n+1}$. Let $f: (\bC^{n+1}, 0) \to (\bC, 0)$ be a non-constant holomorphic function germ such that the Jacobian ideal $J_f$ is contained in $I$ and $j(f) := \dim_\bC I/J_f$ is finite.
    
    Suppose we have germs $G_1, \ldots, G_c, F: (\bC^{n+1+u}, 0) \to (\bC, 0)$ of deformations of $g_1, \ldots, g_c, f$ respectively such that $J_{F \times \pi}$ is contained in $\tilde I := (G_1, \ldots, G_c)$, where $\pi: (\bC^{n+1+u}, 0) \to (\bC^u, 0)$ is the projection. Then the normal cone in $(\bC^{n+1+u}, 0)$ to $(\Sigma_{F \times \pi}, 0)$ is flat over $(\bC^u, 0)$ everywhere on $\pi^{-1}(0) = (\bC^{n+1} \times 0, 0)$, except possibly at the origin if $j(f) > 0$.
\end{prop}

\begin{description}
    \item[Proof] Let $(V(\tilde I), 0) \subset (\bC^{n+1+u}, 0)$ denote the complex-analytic space germ cut out by $\tilde I$. Then, by the Principal Ideal Theorem (e.g., Theorem 10.2 of \cite{eisenbud}), the dimension of $(V(\tilde I), 0)$ is at least $n+1+u-c$; on the other hand, since its fiber over $0 \in \bC^u$ is the vanishing of the ideal $I$ in $(\bC^{n+1}, 0)$ and hence has pure dimension $n+1-c$, we can see by, e.g., Theorem 10.10 of \cite{eisenbud} that its dimension is at most $n+1+u-c$ as well and so it is a complete intersection of codimension $c$ in $(\bC^{n+1+u}, 0)$. In particular, it is Cohen-Macaulay, and hence its flatness over $(\bC^u, 0)$ follows from the fiber dimension by well-known properties of maps from Cohen-Macaulay schemes to smooth bases --- see, e.g., Theorem 18.16 of \cite{eisenbud}.
    
    We now compare $(V(\tilde I), 0)$ and $(\Sigma_{F \times \pi}, 0)$. Consider the following short exact sequence of $\cO_{\bC^{n+1+u},0}$-modules: $$0 \to \tilde I/J_{F \times \pi} \to \cO_{\Sigma_{F \times \pi},0} \to \cO_{V(\tilde I),0} \to 0$$ By the flatness of $(V(\tilde I), 0)$ over $(\bC^u, 0)$ and the long exact sequence in Tor, we can see that this sequence remains exact on restriction over $0 \in \bC^u$. As such, the fiber over $0 \in \bC^u$ of $\tilde I/J_{F \times \pi}$ is precisely $I/J_f$. Our requirement that $j(f) < \infty$ tells us exactly that $I/J_f$ has finite length --- in particular, its support is zero-dimensional, hence concentrated at the origin in $\bC^{n+1} \times 0$, and indeed it is empty if we have $j(f) = 0$. As such, Nakayama's Lemma (e.g., Corollary 4.8 of \cite{eisenbud}) tells us for any sufficiently small representative that the stalks of $\tilde I/J_{F \times \pi}$ will be zero at all points of $\bC^{n+1} \times 0$, except possibly at the origin if $j(f) > 0$. Thus $\Sigma_{F \times \pi}$ agrees with $V(\tilde I)$ at such points and so it will be enough to show that the normal cone to $(V(\tilde I), 0)$ in $(\bC^{n+1+u}, 0)$ is flat over $(\bC^u, 0)$.
    
    Since $(V(\tilde I), 0)$ is a complete intersection in $(\bC^{n+1+u}, 0)$, however, its normal cone in $(\bC^{n+1+u}, 0)$ is a vector bundle by Proposition \ref{prop:cinc} and hence flat over it. As such, the result follows from the flatness of $(V(\tilde I), 0)$ over $(\bC^u, 0)$, which we have already shown.
\end{description}

This result can, of course, be applied with Theorem \ref{thm:strat} to obtain stratifications satisfying the Thom condition with respect to particular pairs of strata, in the spirit of Section 4 of \cite{gaffney_pairs}. It is then immediate that such deformations preserve the smooth fiber --- we state this fact separately for easy reference:

\begin{coro} \label{cor:ci}
    Let $I = (g_1, \ldots, g_c) \subset \cO_{\bC^{n+1}, 0}$ be an ideal defining a germ at the origin of a complete intersection of codimension $c$ in $\bC^{n+1}$. Let $f: (\bC^{n+1}, 0) \to (\bC, 0)$ be a non-constant holomorphic function germ such that the Jacobian ideal $J_f$ is contained in $I$ and $j(f) := \dim_\bC I/J_f$ is finite.
    
    Suppose we have germs $G_1, \ldots, G_c, F: (\bC^{n+1+u}, 0) \to (\bC, 0)$ of deformations of the map-germs $g_1, \ldots, g_c, f$ respectively such that $J_{F \times \pi}$ is contained in $\tilde I := (G_1, \ldots, G_c)$, where $\pi: (\bC^{n+1+u}, 0) \to (\bC^u, 0)$ is the projection. Then the conclusion of Theorem \ref{thm:fibr} holds for $F \times \pi$ --- that is, if we fix a small enough $\varepsilon > 0$, we can choose $\delta > 0$ and $\gamma > 0$ small enough relative to $\varepsilon$ so that the following are true: The representatives $F: B_\varepsilon \times B_\gamma \to \bC$ and $\Sigma_{F \times \pi} \subset B_\varepsilon \times B_\gamma$ exist and, if we let $\Delta := (F \times \pi)(\Sigma_{F \times \pi}) \subset \bC \times B_\gamma$ be the discriminant, the restriction $$F \times \pi: (B_\varepsilon \times B_\gamma) \cap (F \times \pi)^{-1}((D_\delta \times B_\gamma) \setminus \bar \Delta) \to (D_\delta \times B_\gamma) \setminus \bar \Delta$$ defines a smooth locally trivial fibration, where $B_\varepsilon \subset \bC^{n+1}$, $D_\delta \subset \bC$, and $B_\gamma \subset \bC^u$ are the open balls of the specified radii.
\end{coro}

\begin{description}
    \item[Proof] This follows from Proposition \ref{prop:ci} and Condition (ii) of Theorem \ref{thm:fibr}.
\end{description}

Hence a function with scheme-theoretic critical locus a complete intersection can be studied in terms of deformations of that complete intersection, provided compatible deformations of the function can be found. This is easiest in the case where the function, $f$, lies in the square of the ideal $I$ defining the complete intersection, since we can then write $f = \vec g^T H \vec g$ for $\vec g$ a vector with entries the minimal generators of $I$ and simply deform $H$ and $\vec g$ separately to produce such a deformation of $f$. However, the following proposition shows that this in fact imposes fairly stringent restrictions on the complete intersection in question --- compare Theorem 1.13 of \cite{pell_1d}, Theorem 2.2 of \cite{zah_2d}, and Lemma 4.4 of \cite{gaffney_pairs}.

\begin{prop} \label{prop:icis}
    Let $I = (g_1, \ldots, g_c) \subset \cO_{\bC^{n+1}, 0}$ be an ideal defining a germ at the origin of a complete intersection of codimension $c \le n$ in $\bC^{n+1}$. Let $f: (\bC^{n+1}, 0) \to (\bC, 0)$ be a holomorphic function germ such that the Jacobian ideal $J_f$ is contained in $I$ and $j(f) := \dim_\bC I/J_f$ is finite. Then the following are equivalent:
    \begin{enumerate}[label=(\roman*)]
        \item $f \in I^2$.
        \item $I$ is radical.
        \item $I$ defines a reduced isolated complete intersection singularity.
    \end{enumerate}
\end{prop}

\begin{description}
    \item[Proof] (iii) implies (ii) by definition. To see that (ii) implies (i), we use the fact that $f \in \sqrt{J_f}$ (because $0 \in \bC$ is an isolated critical value of $f$ and we work locally) --- since $J_f \subseteq I$, we then have $f \in \sqrt{I} = I$. Thus, by Proposition 1.8 of \cite{pell_findet}, $J_f \subseteq I$ implies that $f \in I^2$.
    
    To show that (i) implies (iii), we suppose that $f \in I^2$ and write $f = \vec g^T H \vec g$ for $$\vec g = \begin{bmatrix} g_1 \\ \vdots \\ g_c \end{bmatrix}$$ and $H$ a $c \times c$ matrix with entries in $\cO_{\bC^{n+1}, 0}$. Then, for any of the differential operators $\partial_i := \tfrac{\partial}{\partial x_i}$, we have
    \begin{alignat*}{2}
        \partial_i f &= \partial_i(\vec g^T H \vec g) \\
        &= (\partial_i \vec g^T) H \vec g + \vec g^T (\partial_i H) \vec g + \vec g^T H (\partial_i \vec g) \\
        &= ((\partial_i \vec g^T) H \vec g)^T + \vec g^T (\partial_i H) \vec g + \vec g^T H (\partial_i \vec g) \\
        &= \vec g^T(H^T(\partial_i \vec g) + (\partial_i H) \vec g + H (\partial_i \vec g)) \\
        &= \vec g^T((H^T + H)(\partial_i \vec g) + (\partial_i H) \vec g).
    \end{alignat*}
    Now, if we let $R := \cO_{\bC^{n+1}, 0}$, we can see that the quotient map $R \to R/J_f$ is the cokernel of the map $R^{\oplus (n+1)} \to R$ given by the Jacobian matrix of $f$, the $i$th column of which is $\partial_i f$. The prior computation then shows that we can decompose this map as $R^{\oplus(n+1)} \xrightarrow{M} R^{\oplus c} \xrightarrow{\vec g^T} R$, where $M$ is the matrix with $i$th column $(H^T + H)(\partial_i \vec g) + (\partial_i H) \vec g$.
    
    Since $j(f) < \infty$, the origin is an isolated point of the support of $I/J_f$. In particular, the map $J_f/IJ_f \to I/I^2$ from the restricted conormal sheaf of the critical locus to the conormal sheaf of the complete intersection germ $(V(I), 0)$ defined by $I$ fails to be an isomorphism only possibly at the origin, since these subspace germs are everywhere else the same. If we regard our map $R^{\oplus c} \xrightarrow{\vec g^T} R$ instead as a map $R^{\oplus c} \to I$, we can see that tensor by $R/I$ gives an isomorphism $(R/I)^{\oplus c} \cong I/I^2$ per Proposition \ref{prop:cinc} since $g_1, \ldots, g_c$ is a regular sequence. It can then be seen that the map $(R/I)^{\oplus(n+1)} \to (R/I)^{\oplus c}$ induced by $M$ factors as $$(R/I)^{\oplus(n+1)} \twoheadrightarrow J_f/IJ_f \to I/I^2 \cong (R/I)^{\oplus c};$$ in particular, its cokernel is supported only possibly at the origin and hence so too the locus in $V(I)$ where the rank of $M$ drops below $c$.
    
    However, if we let $G: (\bC^{n+1}, 0) \to (\bC^c, 0)$ be the map whose components are the $g_i$, we can see from the definition of $M$ that $M$ is equivalent modulo $I$ to the matrix $(H^T + H)\Jac G$, where $\Jac G$ is the Jacobian matrix of $G$. Hence we find that both $H^T + H$ and $\Jac G$ must have full rank everywhere on $V(I)$ except possibly at the origin; by the Jacobian criterion (e.g., Theorem 16.19 of \cite{eisenbud} or Corollary 2 of Section 2.19 of \cite{fischer}), the vanishing of the sum of $I$ and the ideal of maximal minors of $\Jac G$ gives us exactly the locus where $V(I)$ fails to be scheme-theoretically smooth, and so this tells us that $I$ defines an isolated complete intersection singularity which is moreover reduced except possibly at the origin. Since $c \le n$ and so $V(I)$ has dimension at least $1$, the fact that complete intersections lack embedded components tells us that it is reduced at the origin as well. This concludes the proof.
\end{description}

\begin{rmk}
    This result suggests also that, for a complete intersection of positive dimension, appearing scheme-theoretically as the critical locus of a function should actually be a fairly strong condition --- there are certainly, for example, many reduced complete intersections which have more than simply isolated singularities, but the proposition shows that they cannot arise in this context.
\end{rmk}

That the requirement $f \in I^2$ forces $V(I)$ to be an isolated complete intersection singularity should actually in some sense not be too surprising given our general philosophy and results so far --- in this context we can produce deformations of $f$ compatible with {\it any} given deformation of $V(I)$, and Corollary \ref{cor:ci} shows that the local smooth fibers should remain consistent over the parameter space in such deformations. However, if we believe that the structure of the local smooth fibers should reflect that of the critical loci, we should expect the local type of the deformed $V(I)$ to be constant over the parameter space as well in every such deformation --- in particular, there is some heuristic reason to believe that $V(I)$ itself should admit a Milnor fibration, something which is true of isolated complete intersection singularities (see, e.g., \cite{looijenga}) but not of complete intersections in general.

Since we have a programmatic way to produce good deformations of our function in the situation where the conditions of Proposition \ref{prop:icis} hold, it should be possible to compute the homology of the Milnor fiber in general in these cases --- and, indeed, Zaharia has already done so for $c < n$ in Proposition 5.4 of \cite{zah_2d}, using a more direct approach to verify the consistency of the smooth fibers (see also Theorem 3.12 of \cite{gaffney_pairs} for Gaffney's computation of the invariant needed to use Zaharia's formula in practice). Thus the main use of Corollary \ref{cor:ci} will be in moving beyond the case of reduced complete intersections, which has been the focus of much of the existing literature --- as mentioned, however, this will require a more involved approach to obtaining admissible deformations, which we will not pursue further here.

We conclude with a brief discussion of the case of isolated singularities. As mentioned in Section \ref{sec:intro}, every perturbation of a function defining an isolated hypersurface singularity preserves the smooth fiber, and it is not difficult to see that such deformations trivially satisfy the hypotheses of Theorem \ref{thm:fibr}. Less obviously, though, it turns out that in this case things will behave nicely even at the origin:

\begin{rmk} \label{rmk:isol}
    If $F: (\bC^{n+1+u}, 0) \to (\bC, 0)$ is a deformation of a function germ $f$ defining an isolated singularity, then applying Proposition \ref{prop:ci} with $I = J_f$ and $\tilde I = J_{F \times \pi}$ tells us that the normal cone to the critical locus is flat over the parameter space, even at the origin.
\end{rmk}

This fact is not at all necessary in showing that these deformations are well-behaved, but we include it here as a curiosity.

\section{Comparison with Existing Techniques} \label{sec:comp}

As outlined in Section \ref{sec:intro}, Theorem \ref{thm:fibr} is the latest in a long line of results establishing conditions under which some kind of consistency can be obtained for the smooth fibers of germs of families of holomorphic functions. Here we clarify its relationships with these existing concepts and past treatments of the flatness of the critical locus by producing proofs and examples as appropriate.

We begin with finite determinacy. Specifically, we show that any unfolding through an ideal with respect to which our function has finite extended codimension as defined in \cite{bobadilla_fec} satisfies the hypotheses of our theorem:

\begin{prop} \label{prop:fec}
    Let $n \ge 0$ be an integer, $f: (\bC^{n+1}, 0) \to (\bC, 0)$ a nonzero holomorphic map-germ, and $I \subseteq \cO_{\bC^{n+1}, 0}$ an ideal with respect to which $f$ has finite extended codimension. Then, if $F: (\bC^{n+1+u}, 0) \to (\bC, 0)$ is an $I$-unfolding of $f$ and $\pi: (\bC^{n+1+u}, 0) \to (\bC^u, 0)$ is the projection to the parameter space, the origin will be at worst an isolated point of the failure of flatness of the normal cone $C_{\Sigma_{F \times \pi}}\bC^{n+1+u}$ on $\pi^{-1}(0) = (\bC^{n+1} \times 0, 0)$.
\end{prop}

\begin{description}
    \item[Proof] Taking a representative of $f$ with $\tilde I$ a sheaf of ideals on its domain satisfying $\tilde I_0 = I$, we define as in Section 1 of \cite{bobadilla_fec} a coherent sheaf $$\cF := \frac{\tilde I}{\Theta_{\tilde I, e}(f)}$$ whose concentration at the origin witnesses the finite extended codimension condition. Indeed, the dimension over $\bC$ of the stalk of $\cF$ at any point $p$ of our domain gives the extended $\tilde I_p$-codimension of $f$ as a germ at $p$ --- as such, we find that the extended $\tilde I_p$-codimension of $f$ is zero except possibly when $p = 0$.
    
    Theorem 1.12 of \cite{bobadilla_fec} then tells us that, at every $p \ne 0$ in our representative's domain, the 0-parameter unfolding of $f$ is a versal $\tilde I_p$-unfolding. So long as we choose a sufficiently small representative of $F$, we can see that it will give an $\tilde I_p$-unfolding of the germ of $f$ at each such $p \times 0$. Therefore, by the definition of versality, at each point of $\bC^{n+1} \times 0$ in the domain away from the origin, we can in fact write $F \times \pi$ locally up to isomorphism as $f \times \pi$. In particular, we find in a neighborhood of each such point $p \times 0$ that, up to an automorphism of the ambient space, the inclusion $\Sigma_{F \times \pi} \hookrightarrow \bC^{n+1+u}$ is simply the inclusion $\Sigma_f \times \bC^u \hookrightarrow \bC^{n+1+u}$. Its normal cone is then a product over $\bC^u$ as well, hence flat.
\end{description}

Thus Theorem \ref{thm:fibr} generalizes Theorem 2.2 of \cite{bobadilla_fec}. On the other hand, there are many cases where the hypotheses of Theorem \ref{thm:fibr} are satisfied which do not fall under the rubric of finite extended codimension --- see, e.g., Example \ref{ex:wup}.

The relationship with the L\^{e} numbers is much more tenuous:

\begin{rmk}
    The constancy of the L\^{e} numbers in a family is independent of the hypotheses of Theorem \ref{thm:fibr}. It is not difficult to see that there should be many examples where these hypotheses hold but the L\^{e} numbers at the origin do not remain constant --- even an appropriate deformation by translations of an arbitrarily-chosen nontrivial hypersurface singularity will do the trick, but we can also consider any family of isolated singularities with non-constant Milnor number at the origin (see Remark \ref{rmk:isol}) or, again, Example \ref{ex:wup}.
    
    In the other direction, the one-parameter family of Example 9 in \cite{bobadilla_qs} has constant L\^{e} numbers at the origin, as proven in that paper, but one can easily verify that the part of the non-flat locus of its critical scheme lying over $t = 0$ has positive dimension. 
\end{rmk}

Before moving on to the promised examples, we note that the idea of focusing on deformations which make the critical locus consistent over the parameter space in some flatness-related sense seems to have been considered independently at least twice already --- with the important distinction that the prior works both tend to impose such conditions on the {\it reduced} critical locus. In \cite{defthy}, de Jong and van Straten develop theories of admissible deformations in various related contexts, including the notion of an admissible deformation of a holomorphic function germ in their Definition 1.13. In our language, the combination of $F$ and some flat deformation over the same parameter space of a subspace $\Sigma$ of $V(f)$ will define such a deformation precisely when the deformed $\Sigma$ includes into $\Sigma_{F \times \pi} \cap V(F)$ --- although this construction is flexible enough to accommodate a variety of choices of $\Sigma$, including setting $\Sigma = \Sigma_{F \times \pi}$ as in our techniques, the authors note in their Section 1.A that they are primarily interested in dealing with the flatness of {\it reduced} singular loci, rather than using the analytic structure given by the Jacobian ideal. For our present purposes, this is not a fruitful approach --- Example \ref{ex:wupdef} will show that using the reduced structure leads to deformations which need not preserve the nearby smooth fiber, and hence are of limited interest in the study of the Milnor fibration. Additionally, even if one does take $\Sigma$ to be the non-reduced critical locus, the framework as stated and in particular the requirement that the deformation of $\Sigma$ include into $V(F)$ impose limitations on the types of functions and deformations we can consider --- Example \ref{ex:wup}, for instance, will not satisfy this requirement.

Bobadilla's later paper \cite{bobadilla_topoequising} proves that a related set of conditions, some of them specific to the one-dimensional case, guarantees not just the topological consistency of the Milnor fiber but even the topological $R$-equisingularity of the functions in the family. The precise relationship between these conditions --- in particular, the controls imposed on the generic transversal Milnor numbers and Milnor number at the origin --- and the flatness of the normal cone to the non-reduced critical locus seems worthy of further study, but we will not pursue the matter here.

\section{Examples}
\label{sec:ex}

Here we collect various concrete examples pertaining to topics discussed elsewhere in the paper. We begin with an example use of Theorem \ref{thm:fibr}:

\begin{examp} \label{ex:wup}
    Consider the holomorphic function $f: \bC^3 \to \bC$ given by the polynomial equation $f(x, y, z) = x^3 + xy^2z$ and the two-parameter deformation given by $F((x, y, z), s, t) = (x^2 + y^2z - s)(x - t)$. We will use this deformation to compute the Milnor fiber of $f$ at the origin. To verify that the local smooth fiber at the origin varies consistently in this family, we can apply Theorem \ref{thm:fibr} with either Condition (i) or Condition (ii).
    
    For the former, we note that the critical locus of $f$ is set-theoretically the union of the $y$-axis and the $z$-axis. At all points of this set away from the origin, the germ of the vector field $V = y\frac{\partial}{\partial y} - 2z\frac{\partial}{\partial z}$ will have the desired properties; we can see immediately that $Vf = y\frac{\partial f}{\partial y} - 2z\frac{\partial f}{\partial z} = y(2xyz) - 2z(xy^2) = 2xy^2z - 2xy^2z = 0$, so it remains to verify failure of tangency to spheres centered at the origin at critical points of $f$. For points $p = (0, 0, \zeta)$ with $\zeta \ne 0$, we note that the defining equation for the sphere in question is $x\bar x + y\bar y + z\bar z = |\zeta|^2$ and that $(y\frac{\partial}{\partial y} - 2z\frac{\partial}{\partial z})(x\bar x + y\bar y + z\bar z - |\zeta|^2) = y\bar y - 2z\bar z = |y|^2 - 2|z|^2$ evaluated at $p$ is $-2|\zeta|^2 \ne 0$, so the failure of tangency follows. The verification for points on the $y$-axis is similar. Moreover, noting by Remark \ref{rmk:alg} that we can work algebraically, we compute $$\Tor_1^{\bC[s, t]}\left(\frac{\bC[x, y, z, s, t]}{(\tfrac{\partial F}{\partial x}, \tfrac{\partial F}{\partial y}, \tfrac{\partial F}{\partial z})}, \frac{\bC[s, t]}{(s, t)}\right) = 0,$$ either manually or with the aid of a computer algebra system such as {\sc Macaulay2}, {\sc Singular}, or {\sc Oscar}. For example, one can verify that the Tor module vanishes using the following {\sc Macaulay2} code, which computes a free resolution for the second argument and then tensors it by the first before retrieving the Tor module as the first homology of the resulting complex:
    
    \begin{center}
        \begin{minipage}{390px}
            \small \fontencoding{T1}
            \texttt{i1 : R = QQ[s, t];\\
                i2 : S = R[x, y, z]/minors(1, jacobian((x{\textasciicircum}2 + y{\textasciicircum}2*z - s)*(x - t)));\\
                i3 : C = resolution(R{\textasciicircum}1/ideal(s, t));\\
                i4 : (homology(C ** S))\_1 == 0}
        \end{minipage}
    \end{center}
    
     This will output \texttt{true} if the Tor module is zero, and \texttt{false} otherwise. Note that here we are working over $\bQ$ and using the faithful flatness of $\bC$ over $\bQ$ to note that the corresponding results will hold over $\bC$. Then, by the local criterion (e.g., Theorem 6.8 of \cite{eisenbud}) and the well-behavedness of Tor under the various localizations involved, this is sufficient to guarantee that in particular we have flatness at the origin.
    
    For the latter, we note by Remarks \ref{rmk:alg} and \ref{rmk:algc} that we can again work entirely algebraically, and directly compute $$\Tor_1^{\bC[s, t]}\left(\gr_{\left(\tfrac{\partial F}{\partial x}, \tfrac{\partial F}{\partial y}, \tfrac{\partial F}{\partial z}\right)}\bC[x, y, z, s, t], \frac{\bC[s, t]}{(s, t)}\right) = 0.$$ As before, this is enough to give us the flatness condition we seek.
    
    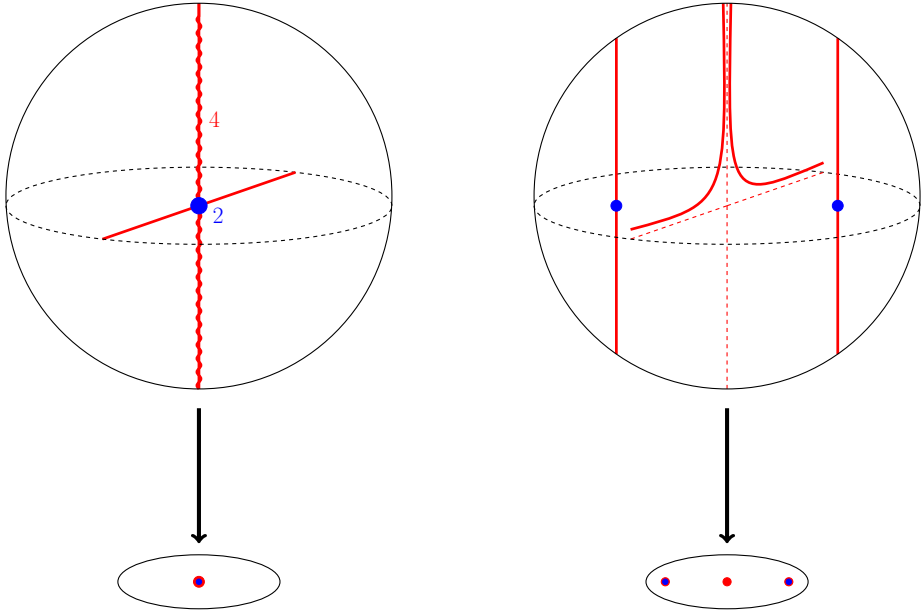
\begin{figure}
        \begin{center}
            \begin{minipage}{150px}
                \resizebox{150px}{!}{\input{Figure_zero}}
            \end{minipage}
            \quad\quad\quad\quad
            \begin{minipage}{150px}
                \resizebox{150px}{!}{\input{Figure_nonzero}}
            \end{minipage}
            \vspace{.1cm}
            \caption{The critical loci of Example \ref{ex:wup} over $t = 0$ and $t \ne 0$ respectively ($s = 5t^2$)}
            \label{fig:wup}
        \end{center}
    \end{figure}
    
    Hence, whichever approach we choose to take, we can conclude by Theorem \ref{thm:fibr} that this deformation preserves the local smooth fiber at the origin. To compute the homology of this fiber, we restrict over the curve $s = 5t^2$ in the parameter space, obtaining a one-parameter deformation of $f$ which we will henceforth call $F$ in a slight abuse of notation. For each $t$, denote by $f_t$ the restriction of $F$ to $\bC^3 \times \{t\}$.
    
    Before we proceed to the computation, since we are interested in having a picture of the scheme-theoretic critical locus and how the smooth fiber can be gleaned from it, we will say a little about how to visualize the critical loci of the functions in this family, as depicted in Figure \ref{fig:wup}. In the $t = 0$ case, we can observe $J_f = (3x^2 + y^2z, 2xyz, xy^2)$ and find by computation that this ideal has associated primes $(x, y)$, $(x, z)$, and $(x, y, z)$ --- it is then straightforward to show that the multiplicities (in the sense of Section 6 of Chapter 3 of \cite{eisenbud}) of these primes in $\bC[x, y, z]/J_f$ are 4, 1, and 2 respectively.
    
    To visualize the $t \ne 0$ case, on the other hand, we find that the ideal $J_{F \times \pi} = (3x^2 + y^2z - 5t^2 - 2tx, 2(x-t)yz, (x-t)y^2)$ cutting out the family of critical loci has associated primes $(x - t, y^2z - 4x^2)$, $(x + t, y)$, $(x + t, y, z)$, $(3x - 5t, y)$, and $(3x - 5t, y, z)$, each corresponding to a component of multiplicity 1. Now, since the isomorphism $\bC^3 \times \bC^* \cong \bC^3 \times \bC^*$ which sends $(x, y, z, t)$ to $(tx, ty, z, t)$ gives a trivialization of the vanishing of $J_{F \times \pi}$ in $\bC^3 \times \bC^*$ over $\bC^*$, we find for any fixed $t \ne 0$ that the components and multiplicities of $\Sigma_{f_t}$ are given simply by restricting the aforementioned components over this value. In particular, each of the {\it minimal} primes in this list corresponds to a connected component of $J_{f_t}$, and these components are sent by $f_t$ to $0$, $8t^3$, and $-\tfrac{40}{27}t^3$ respectively.
    
    Now let $\varepsilon > 0$ be arbitrary and consider the restriction of $F$ to $B_\varepsilon \times \bC$, for $B_\varepsilon$ the open ball of radius $\varepsilon$ centered at the origin in $\bC^3$. Then, provided that $t$ is small enough relative to $\varepsilon$, the critical locus of $f_t$ in $B_\varepsilon$ will be as depicted in Figure \ref{fig:wup}, since the formation of the critical locus commutes with restriction to open sets.
    
    Although the deformation we have used is not exactly of the type considered in \cite{st_deform}, we can still use a slight modification of the methods of that paper, which were originally introduced in \cite{siersma_icis}, to compute the homology of the smooth fiber. Specifically, we observe that, if we take the inverse image under $f$ of a small enough disk $D_\delta$ around the origin and intersect it with $B_\varepsilon$, the result will be homotopy-equivalent to the fiber $f^{-1}(0) \cap B_\varepsilon$, hence contractible. By including the boundary and applying Lemma \ref{lem:tfil}, we find that, for small enough $t \ne 0$, the space $\bE := {f_t}^{-1}(D_\delta) \cap B_\varepsilon$ is contractible as well. Letting $\bF$ be a smooth fiber of $f_t$ in $B_\varepsilon$, we then find that $H_{*+1}(\bE, \bF) \cong \tilde H_*(\bF)$ by the long exact sequence in homology, and Theorem \ref{thm:fibr} tells us that $\tilde H_*(\bF)$ is precisely the reduced homology of the Milnor fiber of $f$.
    
    As in \cite{st_deform}, we can now decompose this relative homology group into a direct sum of the analogous groups around each critical value individually using homotopy retraction and excision. Specifically, we proceed as follows. Let $D_0, D_+, D_-$ be sufficiently small disjoint open disks inside $D_\delta$ around the critical values  $0$, $8t^3$, and $-\tfrac{40}{27}t^3$ respectively. Observing that the transversal Milnor fibration given by fixing $z$ is a local product at all points of the critical locus where $z \ne 0$, we find that the special points (in the sense of \cite{st_deform}) are exactly $q_+ := (-t, 0, 0)$ and $q_- := (\tfrac53t, 0, 0)$; let $B_+$ and $B_-$ be sufficiently small Milnor balls in $B_\varepsilon$ around these points respectively. Also let $\cU_0, \cU_+, \cU_-$ be sufficiently small tubular neighborhoods in $B_\varepsilon$ of the connected components $\Sigma_0 := V(x - t, y^2z - 4x^2) \cap B_\varepsilon, \Sigma_+ := V(x + t, y) \cap B_\varepsilon, \Sigma_- := V(3x - 5t, y) \cap B_\varepsilon$ of the critical locus. Then, if we let $\bE_0 := {f_t}^{-1}(D_0) \cap \cU_0$, $\bE_+ := {f_t}^{-1}(D_+) \cap (\cU_+ \cup B_+)$, and  $\bE_- := {f_t}^{-1}(D_-) \cap (\cU_- \cup B_-)$ and $\bF_0, \bF_+, \bF_-$ be the intersections of the smooth fibers of $f_t$ over small non-critical values nearby the appropriate critical values with $\bE_0, \bE_+, \bE_-$ respectively, we have $$H_{*+1}(\bE, \bF) \cong H_{*+1}(\bE_0, \bF_0) \oplus H_{*+1}(\bE_+, \bF_+) \oplus H_{*+1}(\bE_-, \bF_-)$$ by homotopy retraction and excision.
    
    We can now consider each of these relative homology groups independently. For each of the groups $H_{*+1}(\bE_\pm, \bF_\pm)$, we can take a further retraction of pairs to get the relative homology of ${f_t}^{-1}(D_\pm) \cap B_\pm$ and its intersection with a nearby smooth fiber; by our reasoning above, this is simply the reduced homology of the Milnor fiber of $f_t$ at the point $q_\pm$. Since $x - t$ is a unit in the local power series ring at this point and multiplication by units does not affect the Milnor fiber, the point is a $D_\infty$ singularity and hence the reduced homology of the Milnor fiber consists of a single $\bZ$-summand in degree 2.
    
    Consider the remaining component $\Sigma_0$ of the critical locus and note that $\bE_0 \simeq \Sigma_0 \simeq S^1$. Thus, from the long exact sequence of a pair in homology and the fact that the bounds of \cite{katmat73} imply $H_{i+1}(\bE, \bF) \cong \tilde H_*(\bF) \cong 0$ and hence $H_{i+1}(\bE_0, \bF_0) \cong 0$ for $i \not\in \{1, 2\}$, we obtain a short exact sequence $$0 \to H_2(\bE_0, \bF_0) \to H_1(\bF_0) \to H_1(\bE_0) \to 0$$ with $H_1(\bE_0) \cong \bZ$ and an isomorphism $H_3(\bE_0, \bF_0) \cong H_2(\bF_0)$.
    
    Now, by considering the transversal slices given by the level sets of the linear function $y$, we see that the transversal type of the singularity is everywhere $A_1$ on $\Sigma_0$ and the transversal Milnor fibration is a local product as before. Indeed, by noting $y \ne 0$ around $\Sigma_0$ and writing $f_t = (\tilde x^2 + \tilde z - 5t^2)(\tilde x - t)$ in the coordinates $(\tilde x, \tilde y, \tilde z) = (x, y, y^2z)$, we obtain a global trivialization of the transversal Milnor fibration, so $\bF_0$ is homotopy equivalent to a torus and the map $S^1 \times S^1 \simeq \bF_0 \to \bE_0 \simeq S^1$ is simply the projection onto one factor --- here we are able to ignore any issues near the boundary of $B_\varepsilon$ by an appropriate homotopy retraction.
    
    As such, $H_3(\bE_0, \bF_0) \cong H_2(\bF_0) \cong \bZ$, and the map $\bZ \oplus \bZ \cong H_1(\bF_0) \to H_1(\bE_0) \cong \bZ$ is the projection onto a single factor and hence has kernel $H_2(\bE_0, \bF_0) \cong \bZ$. Combining this with our previous computation, we find that the reduced homology of the Milnor fiber of $f$ is given by $$\tilde H_i(\bF) \cong \begin{cases}
        \bZ & i = 1 \\
        \bZ^{\oplus 3} & i = 2 \\
        0 & \text{else.}
    \end{cases}$$
\end{examp}

Thus we can apply Theorem \ref{thm:fibr} to compute Milnor fiber homology by splitting off not only Morse points but also higher-dimensional features, such as the two $D_\infty$ singularities seen above.

Now consider Conditions (i) and (ii) of Theorem \ref{thm:fibr}. The condition of flatness of the critical locus in a family, as in Condition (i), is easier to verify than the condition of the flatness of its normal cone as in Condition (ii) --- as discussed in Section \ref{sec:intro}, the difference between these cases is essentially whether we require consistency of the critical locus merely in an intrinsic sense or also in an embedded one. The idea that this intrinsic consistency might be sufficient is theoretically appealing, and moreover would greatly reduce the computational cost of these methods in practice --- in the procedure of Corollary \ref{cor:homstrat}, for example, the computation of the normal cone can be quite expensive. However, the following example shows that not every function satisfies the hypothesis on the existence of vector field germs necessary for Condition (i) to hold:

\begin{examp} \label{ex:novec}
    Let $f: (\bC^3, 0) \to (\bC, 0)$ be the holomorphic function germ given by $f(x, y, z) = x^6 + y^6 + x^2y^2z^2$. Then we claim that no deformation of $f$ can satisfy Condition (i) of Theorem \ref{thm:strat}.
    
    To show this, observe that the critical locus of $f$ is set-theoretically the $z$-axis in $\bC^3$ and consider a point $p = (0, 0, \zeta)$ for $\zeta \ne 0$. Suppose toward a contradiction that there exists a vector field germ $V$ at $p$ satisfying $Vf = 0$ such that $V(p)$ is not contained in the tangent space to $S_{|p|}$ at $p$. Then we see that $(dz \cdot V)(p) \ne 0$, so we have $V = a \frac{\partial}{\partial x} + b \frac{\partial}{\partial y} + u \frac{\partial}{\partial z}$ for holomorphic function germs $a, b, u \in \cO_{\bC^3, p}$ with $u$ a unit. In particular, since $Vf = 0$, we can see that $\frac{\partial f}{\partial z}$ is contained in the ideal of this convergent power series ring generated by $\frac{\partial f}{\partial x}$ and $\frac{\partial f}{\partial y}$.
    
    Letting $J_f$ be the Jacobian ideal of $f$ in $\bC[x, y, z]$, $I$ the ideal generated by $\frac{\partial f}{\partial x} = 6x^5 + 2xy^2z^2$ and $\frac{\partial f}{\partial y} = 6y^5 + 2x^2yz^2$, and $\mFm = (x, y, z - \zeta)$, we can then see by the faithful flatness of the map $\bC[x, y, z]_\mFm \to \cO_{\bC^3, p}$ (as in Remark \ref{rmk:alg}) that $(J_f/I)_\mFm$ must be zero. However, a straightforward computation in a computer algebra system such as {\sc Macaulay2}, {\sc Singular}, or {\sc Oscar} reveals that this is not the case, yielding the desired contradiction.
    
    However, Theorem \ref{thm:fibr} can still be applied to deformations of $f$ using Condition (ii). For example, consider the one-parameter deformation of $f$ given by $F((x, y, z), t) = x^6 + y^6 + x^2y^2(z^2 - t)$. By again using a computer algebra system, we can easily verify that $t$ is not a zerodivisor on $\gr_{J_{F \times \pi}} \bC\{x, y, z, t\}$ (with $\pi$ the projection as usual) and hence that $C_{\Sigma_{F \times \pi}} \bC^4$ is flat over the parameter space --- as such, Theorem \ref{thm:fibr} tells us that $F$ can be used to study the Milnor fibration of $f$. Although this deformation does not simplify the critical locus far enough to make a complete computation of the Milnor fiber's homology tractable without significant further work, one can see by computing a primary decomposition of $J_{F \times \pi}$ that it splits off 12 Morse points, given by the primary ideals $(x - y, 3y^2 - t, z)$, $(x + y, 3y^2 - t, z)$, $(3xy - t, 3x^2 + 3y^2 + t, z)$, and $(3xy + t, 3x^2 + 3y^2 + t, z)$; therefore, the methods of \cite{st_deform}, used as in Example \ref{ex:wup}, tell us immediately that the rank of $\tilde H_2(\bF_f)$ is at least 12.
\end{examp}

From this example we can see that our current proof methods are not enough to establish the sufficiency of the intrinsic condition in general, and it seems ambiguous at this point whether we should expect it to be true. One natural approach to attempting to prove or disprove it would of course be to search for deformations of functions such as the one in the previous example which have flat critical locus but not flat normal cone --- however, in practice even finding interesting deformations of a given function simply satisfying the intrinsic condition can be quite difficult, so we defer such investigations to a later work.

As promised in Section \ref{sec:comp}, we now give a family where the {\it reduced} critical locus is deformed flatly, as suggested in \cite{defthy} (although that paper is primarily interested in cases where the space of such ``admissible deformations'' is finite-dimensional, which we do not claim here), as is its normal cone, but the smooth fibers are not consistent even in very coarse senses:

\begin{examp} \label{ex:wupdef}
    Let $f: (\bC^3, 0) \to (\bC, 0)$ be the holomorphic function germ given by $f = y^2 - x^3$ and consider the one-parameter deformation of $f$ given by $F = y^2 + x^2(tz - x)$. Then, in the language of Theorem \ref{thm:fibr}, $\Sigma_{F \times \pi}$ has a non-isolated failure of flatness over the origin in the parameter space. However, we can see that the trivial deformation of the reduced critical locus $V(x, y)$ of $f$ does include into $\Sigma_{F \times \pi}$, since $F \in (x, y)^2$, so $(V(x, y) \times \bC, F)$ is an admissible deformation over the parameter space $\bC$ in the sense of Definition 1.13 of \cite{defthy} when the reduced critical locus is used. Indeed, even the normal cone to the subspace cut out by $\sqrt{J_{F \times \pi}}$ is a product over the parameter space, hence flat.
    
    Now, the Milnor fiber of $f$ is given up to homotopy by the Milnor fiber for the germ $(\bC^2, 0) \to (\bC, 0)$ given by the same expression $y^2 - x^3$, since if we are looking at the restriction of $f$ to $B_\varepsilon \cap f^{-1}(D_\delta^*)$ for $\varepsilon, \delta > 0$ small enough to define the Milnor fibration we can simply collapse the vertical coordinate to get the 2-dimensional picture. Thus, since the Milnor number of this cusp is 2, we find that the Milnor fiber of $f$ has the homotopy type of the bouquet $S^1 \vee S^1$.
    
    On the other hand, keeping the same $\varepsilon$, we can see for any $t \ne 0$ that the restriction $f_t$ of $F$ to $B_\varepsilon \times \{t\}$ gives a $D_\infty$ singularity at the origin --- the isomorphism to the Whitney umbrella can be obtained by a straightforward change of coordinates. By well-known properties of the Whitney umbrella, we find immediately that the stratification given by the filtration $\bC^3 \supset V(f_t) \supset V(x, y) \supset V(x, y, z)$ satisfies the Thom ($a_{f_t}$) condition; thus, as mentioned in Remark \ref{rmk:radius}, we can compute the Milnor fibration of $f_t$ using $B_\varepsilon$ so long as $S_\varepsilon$ and all smaller spheres at the origin are transverse to its strata. This is immediate for all strata except the regular part of $V(f_t)$.
    
    On this stratum, failure of transversality will occur precisely when the $3 \times 6$ real matrix whose rows are the gradients of the real and imaginary parts of $f_t$ and the squared distance from the origin fails to be full rank. By reinterpreting every second column as an imaginary coordinate and using the Cauchy-Riemann equations, we can see that, away from the origin, this condition is equivalent to the requirement that the vector of complex conjugates of partial derivatives of $f_t$ --- the gradient of $f_t$, in the language of \cite{milnor68} --- is a complex multiple of the vector $(x, y, z)$ of coordinates of the ambient space. Some calculations, which we omit for the sake of brevity, then tell us that this failure of transversality can be excluded by keeping $t$ small enough relative to $\varepsilon$.
    
    Hence the smooth fiber of $f_t$ in $B_\varepsilon \times \{t\}$ over any value near enough to the origin in $D_\delta$ will be the Milnor fiber of a $D_\infty$ singularity and so have the homotopy type of $S^2$. This is not consistent with the Milnor fiber at $t = 0$ even at the level of Euler characteristic, much less homology, homotopy type, or diffeomorphism type.
\end{examp}

Thus deformations of this sort are of limited utility in obtaining Milnor fiber information --- for such purposes, it really is necessary to keep track of the non-reduced structure of the critical locus, as in Theorems \ref{thm:strat} and \ref{thm:fibr}.

\section{Bibliography}

\printbibliography[heading=none]
\end{document}

%% file: Figure_zero.tex
\begin{tikzpicture}
    \draw [thick, dashed] (5, -.25) arc [start angle=0, end angle=180, x radius=5cm, y radius=1cm];
    
    \draw [line width=2pt, color=red] (0, -.25) -- +(60:5cm and 1cm);
    \draw [line width=2pt, color=red] (0, -.25) -- +(240:5cm and 1cm);
    
    \draw [line width=2pt, color=red, decorate, decoration={coil,aspect=0,amplitude=1.5pt}] (0, -5) -- (0, 5);
    \draw [line width=2pt, color=red] (0, -5) -- (0, 5);
    \draw [color=red] (.4, 2) node {\Large 4};
    
    \filldraw [color=blue] (0, -.25) circle (6pt);
    \draw [color=blue] (.5, -.5) node {\Large 2};
    
    \draw [thick] (5, 0) arc [start angle=0, end angle=360, radius=5cm];
    \draw [thick, dashed] (-5, -.25) arc [start angle=180, end angle=360, x radius=5cm, y radius=1cm];
    
    \draw [thick] (2.1, -10) arc [start angle=0, end angle=360, x radius=2.1cm, y radius=.7cm];
    
    \draw[->, line width=3pt] (0, -5.5) -- (0, -9);
    
    \filldraw [color=red] (0, -10) circle (4pt);
    \filldraw [color=blue] (0, -10) circle (2pt);
\end{tikzpicture}

%% file: Figure_nonzero.tex
\begin{tikzpicture}
    \draw [thick, dashed] (5, -.25) arc [start angle=0, end angle=180, x radius=5cm, y radius=1cm];
    
    \draw [thick, dashed, color=red] (0, -.25) -- +(60:5cm and 1cm);
    \draw [thick, dashed, color=red] (0, -.25) -- +(240:5cm and 1cm);
    
    \draw [line width=2pt, color=red] (240:5cm and 1cm) .. controls (0, -.2) .. (-.1, 5);
    \draw [line width=2pt, color=red] (.1, 5) .. controls (0, -.2) .. (60:5cm and 1cm);
    
    \draw [thick, dashed, color=red] (0, -5) -- (0, 5);
    
    \draw [line width=2pt, color=red] (125:5cm) -- (235:5cm);
    \filldraw [color=blue] (-2.868, -.25) circle (4pt);
    
    \draw [line width=2pt, color=red] (55:5cm) -- (305:5cm);
    \filldraw [color=blue] (2.868, -.25) circle (4pt);
    
    \draw [thick] (5, 0) arc [start angle=0, end angle=360, radius=5cm];
    \draw [thick, dashed] (-5, -.25) arc [start angle=180, end angle=360, x radius=5cm, y radius=1cm];
    
    \draw [thick] (2.1, -10) arc [start angle=0, end angle=360, x radius=2.1cm, y radius=.7cm];
    
    \draw[->, line width=3pt] (0, -5.5) -- (0, -9);
    
    \filldraw [color=red] (0, -10) circle (3pt);
    
    \filldraw [color=red] (-1.6, -10) circle (3pt);
    \filldraw [color=blue] (-1.6, -10) circle (2pt);
    
    \filldraw [color=red] (1.6, -10) circle (3pt);
    \filldraw [color=blue] (1.6, -10) circle (2pt);
\end{tikzpicture}